\documentclass[preprint,12pt]{elsarticle}

\usepackage{amssymb,amsmath,bm}
\usepackage{hyperref}

\biboptions{numbers,sort&compress}  

\journal{Computer Methods in Applied Mechanics and Engineering}

\begin{document}

\begin{frontmatter}

\title{Preconditioning Techniques for Hybridizable Discontinuous Galerkin Discretizations on GPU Architectures}

\author[mit]{Andrew Welter}
\author[mit]{Ngoc Cuong Nguyen}
\address[mit]{Department of Aeronautics and Astronautics, Massachusetts Institute of Technology, \ 77 Massachusetts Avenue, Cambridge, MA 02139, USA}

\begin{abstract}
We present scalable iterative solvers and preconditioning strategies for Hybridizable Discontinuous Galerkin (HDG) discretizations of partial differential equations (PDEs) on graphics processing units (GPUs). The HDG method is implemented using GPU-tailored algorithms in which local element degrees of freedom are eliminated in parallel, and the globally condensed system is assembled directly on the device using dense-block operations. The global matrix is stored in a block format that reflects the natural HDG structure, enabling all iterative solver kernels to be executed with strided batched dense matrix–vector multiplications. This implementation avoids  sparse data structures, increases arithmetic intensity, and sustains high memory throughput across a range of meshes and polynomial orders. The nonlinear solver combines Newton’s method with preconditioned GMRES, integrating scalable preconditioners such as block-Jacobi, additive Schwarz domain decomposition, and polynomial smoothers. All preconditioners are implemented in batched form with architecture-aware optimizations—including dense linear algebra kernels, memory-coalesced vector operations, and shared-memory acceleration—to minimize memory traffic and maximize parallel occupancy. Comprehensive studies are conducted for a variety of PDEs (including Poisson equation, Burgers equation, linear and nonlinear elasticity, Euler equations, Navier-Stokes equations, and Reynolds–Averaged Navier–Stokes equations) using structured and unstructured meshes with different element types and polynomial orders on both NVIDIA and AMD GPU architectures. 
\end{abstract}

\end{frontmatter}


\section{Introduction}

Hybridizable Discontinuous Galerkin (HDG) methods have emerged as high-order finite element discretizations for the numerical solution of partial differential equations (PDEs) in a wide variety of applications \cite{Nguyen2012,Schutz2013,Woopen2014c,Fidkowski2016,Williams2017, Nguyen2023c,Ciuca2020,Vidal-Codina2021,Sanchez2021,Solano2022,StanglmeierNguyenPeraireCockburn16,Terrana2019a,Nguyen2011j,Huynh2013c,Nguyen2024d,NguyenPeraireCockburnHDGNS11,Ueckermann2016,SANCHEZ2022114969}. By introducing additional degrees of freedom on the mesh skeleton, HDG methods enable local static condensation and reduce the global linear system to the trace unknowns \cite{CockburnGopalakrishnanLazarov09,NguyenPeraireCockburnConvDiff09,nguyen2009,Nguyen2011a,Cockburn2010e,Nguyen2010}. This property not only improves  accuracy and efficiency, but also lends itself naturally to parallel implementation \cite{VilaPerez2022,Huerta2013,N2013,Nguyen2023a,Hoskin2024,CuongNguyen2022, Nguyen2020gpu}. Despite these advantages, the global linear systems resulting from HDG discretizations can be large and expensive to solve for large-scale problems. Direct solvers are often impractical due to memory and time constraints, especially in three dimensions with high-order discretizations. Iterative solvers, paired with effective preconditioners, are therefore essential to make HDG methods viable for large-scale simulations \cite{Fernandez2017a,CuongNguyen2022,Hoskin2024}.

There is a rich body of work on iterative solvers and preconditioning techniques for DG and HDG methods. Popular choices for preconditioners include block-Jacobi, block Gauss-Seidel, and block-ILU factorization \citep{persson:GMRESDG}. Early efforts focused on block-Jacobi and Gauss-Seidel smoothers in multigrid frameworks for DG discretizations of linear problems \citep{Gopalakrishnan2003,Nastase2006}. For nonlinear problems, including the compressible Navier–Stokes equations, block tridiagonal and block-ILU smoothers have been shown to improve convergence and robustness, particularly in the presence of grid stretching and high Reynolds numbers \citep{Fidkowski2005,persson:GMRESDG}. Although block-ILU is more difficult to implement, it is significantly more effective. Its performance depends critically on the ordering of unknowns, for which minimum discarded fill algorithms are often employed. Multigrid methods, both $h$- and $p$-multigrid, have also been developed to solve linear systems from DG/HDG discretizations \citep{Gopalakrishnan2003,Nastase2006,Fidkowski2005,Shahbazi2009,Franciolini2020,Bastian2019,Betteridge2021,Wildey2019}.  Fernandez et al. \citep{Fernandez2017a} developed a preconditioned Newton–GMRES solver for HDG discretizations, combining a restricted additive Schwarz (RAS) method with block-ILU(0) smoothing and a graph-based reordering strategy. This approach enabled efficient large eddy simulations of transitional flows at high Reynolds numbers. Additional strategies such as algebraic multigrid (AMG), Schur complement reduction, deflation, and polynomial preconditioners have been explored to improve solver scalability for HDG linear systems \citep{Muralikrishnan2017}.  More recently, Nguyen et. al. \citep{CuongNguyen2022, VilaPerez2022} proposed a matrix-free reduced basis preconditioner for JFNK solvers, enabling the simulation of transonic buffet phenomena in high Reynolds regimes using implicit DG solvers. 

However, since many of the existing implicit DG solvers are implemented on CPU thread models and hierarchical memory layouts, they  exhibit poor performance on GPUs due to synchronization overheads, irregular memory access, and limited concurrency. As modern scientific computing increasingly relies on GPU-based hardware, the need for iterative solvers and preconditioners that are  effective on GPUs has become critical. Classical techniques such as ILU and multigrid involve sparse triangular solves and coarse-grid corrections that are difficult to parallelize on GPU architectures. These limitations have spurred interest in GPU-aware solvers and preconditioners.


GPU acceleration has been extensively explored for high-order explicit discontinuous Galerkin (DG) methods. Early work by Klöckner et al.\ \citep{klockner:DGGPU} demonstrated the suitability of DG discretizations for GPUs, achieving high performance via element-local operations and memory-aware kernel design. Subsequent studies extended GPU-accelerated DG to a broad range of applications, including seismic wave propagation on unstructured tetrahedral meshes \citep{Mu2013}, nonlinear shallow water equations with positivity-preserving limiters and local time-stepping \citep{Gandham2015}, compressible Navier--Stokes simulations on hybrid unstructured grids using OpenACC directives \citep{Xia2015}, and non-hydrostatic atmospheric dynamics in a  continuous/discontinuous Galerkin framework  through OCCA-based kernels \citep{Abdi2019}. Performance-optimized implementations have been developed for hybrid meshes \citep{chan2016}, Bernstein--Bézier formulations for wave propagation \citep{chan2017}, acoustic and elastic wave models \citep{Modave2016}. More recent efforts focus on GPU-accelerated DG schemes on polytopic meshes \citep{Dong2021}, GPU-based high-order DG solvers for supersonic turbulence in astrophysical flows \citep{Cernetic}, and p-adaptive DG formulations for shallow water equations on heterogeneous CPU–GPU architectures \citep{Faghih2025}. Collectively, these works establish a mature ecosystem of GPU-accelerated technologies and motivate analogous developments for impicit DG solvers.

GPU implementations of implicit DG solvers remain relatively limited compared to their explicit counterparts. Karakus et al.\ \citep{Karakus2019} developed an incompressible Navier--Stokes solver featuring fully GPU-accelerated multigrid preconditioners, demonstrating the feasibility of high-order implicit formulations on modern GPU architectures. Complementary advances in operator efficiency have focused on fast matrix-free DG implementations leveraging sum-factorization, vectorization, and cache-aware data layouts \citep{Kronbichler2019}. For hybridizable DG methods, Roca et al.\ \citep{Roca2011} introduced a custom sparse matrix–vector product tailored to the HDG trace system, outperforming vendor-optimized GPU libraries such as cuSPARSE. Fabien \citep{Fabien2020} developed a GPU-accelerated HDG solver for linear elasticity, using a time–accuracy–size metric and roofline analysis to guide kernel design and performance tuning. Recent work has expanded the scope of implicit DG methods on GPUs to large-scale, high-speed flow simulations. Matrix-free GPU solvers have been developed for transonic buffet \citep{CuongNguyen2022}, hypersonic LES \citep{Nguyen2020gpu, Nguyen2023a}, turbulent and shock-dominated compressible flows \citep{Hoskin2024}, and parametric space-weather and thermospheric modeling \citep{VilaPerez2024, VilaPerez2022}. These studies demonstrate that fully implicit or Jacobian-free Newton--Krylov DG solvers can attain high accuracy and excellent scalability across thousands of GPUs, even for highly stiff, multiscale flow regimes. These contributions provide a foundation for scalable implicit DG and HDG solvers on GPU platforms. However, they also underscore the need for GPU-based preconditioners and iterative strategies for high-order HDG discretizations.


In this work, we develop and implement iterative solvers and preconditioners for HDG discretizations, targeting both NVIDIA and AMD architectures. Our approach integrates Newton-Krylov methods with a suite of preconditioning strategies, including block-Jacobi, polynomial preconditioners, and additive Schwarz methods. These preconditioners are designed for compatibility with GPU hardware, emphasizing memory efficiency, data locality, and high occupancy. We rely primarily on cuBLAS (and hipBLAS) dense linear algebra libraries for all major computational kernels, including local static condensation, local matrix inversion, matrix–matrix and matrix–vector operations, and batched preconditioner applications. Furthermore, we employ Kokkos \cite{Carter2014,Trott2022} to implement all remaining GPU kernels that are not served by vendor-optimized BLAS routines. This includes element-level flux evaluations, Jacobian and residual assembly, trace and volume operator applications, vector updates, and communication packing/unpacking. By relying on Kokkos for these components, we achieve a portable, single-source implementation that runs efficiently on both CUDA and HIP backends.

As a result, much of the lower-level performance tuning—such as thread-block scheduling, memory tiling, and kernel fusion—is handled automatically by the  BLAS libraries and Kokkos backends. This design yields robust and portable performance across a broad range of polynomial orders, element types, and mesh configurations while minimizing the need for architecture-specific optimizations. The combination of cuBLAS/hipBLAS for dense kernels and Kokkos for all remaining GPU-parallel operations allows the HDG solver to exploit the GPU memory hierarchy and arithmetic throughput, while maintaining concise, maintainable, and hardware-agnostic code. 

We evaluate the proposed solvers on a broad and diverse set of PDE models, including the heat equation, convection–diffusion, wave propagation, linear and nonlinear elasticity, compressible Euler and Navier–Stokes equations, and Reynolds-averaged Navier–Stokes (RANS) models. These test cases encompass a wide range of physical regimes—from smooth flows to problems with shocks and boundary layers—as well as structured and unstructured meshes and various polynomial orders. We assess solver performance across  GPU platforms in terms of time-to-solution, iteration counts,  polynomial degrees, mesh sizes, and other parameters.

The contributions of this work are threefold. First, we develop a set of preconditioning strategies that are specifically designed for HDG discretizations on modern GPU architectures. Second, we present a suite of GPU-aware solver implementations that exploit the dense-block structure and face-based coupling inherent in HDG global systems to enable high arithmetic intensity and efficient use of device memory. Third, we demonstrate the performance of the proposed HDG solver across multiple classes of PDEs and hardware backends, including both NVIDIA and AMD GPUs. Although the focus is on HDG methods, the underlying algorithms—particularly the GPU-optimized preconditioners, matrix assembly routines, and Krylov subspace solvers—are broadly applicable to other DG and finite element methods. 

The remainder of this paper is organized as follows. In Section 2, we review the HDG formulation and the structure of the resulting linear systems. Section 3 discusses the iterative solvers and preconditioning techniques considered in this work. Section 4 presents our GPU implementation strategies. Numerical results are given in Section 5. In Section 6 we conclude with a summary and directions for future work.

\section{HDG Discretization and Linear System Structure}

\subsection{HDG Discretization for Systems of Conservation Laws}

We consider a system of conservation laws for $M$ state variables defined on a physical domain $\Omega \subset \mathbb{R}^D$. The time-dependent system of conservation laws is given by
\begin{equation}
\label{eq:conservation}
\frac{\partial \bm{u}}{\partial t} + \nabla \cdot \bm{F}(\bm{u}, \nabla \bm{u}) = \bm s(\bm u, \nabla \bm u) \quad \text{in } \Omega \times (0, T],
\end{equation}
where $\bm u(\bm x, t) \in \mathbb{R}^M$ is the solution of the system of conservation laws,  the  physical fluxes $\bm F = (\bm f_1(\bm u, \nabla \bm u), \ldots, \bm f_D(\bm u, \nabla \bm u))$, of dimension $M \times D$ include $D$ vector-valued functions of the solution and its gradient, and $\bm s(\bm u, \nabla \bm u) \in \mathbb{R}^M$ is the source term. The system is supplemented with the prescription of initial and boundary conditions
\begin{equation}
\begin{split}
\bm{u}(\bm{x}, 0) &= \bm{u}_0(\bm{x}) \quad \text{in }  \Omega, \\ 
\bm{b}(\bm{u}, \nabla \bm{u}, \bm n) & = 0 \qquad \quad \mbox{on }  \partial \Omega \times (0, T],
\end{split}
\end{equation}
where $\bm{u}_0$ denotes an initial state,  $\bm b$ denotes a boundary condition operator, and $\bm n$ is the normal vector at the domain boundary. In the steady-state regime, the governing equations reduce to
\begin{equation}
\label{eqst}
 \nabla \cdot \bm F(\bm u, \nabla \bm u) = \bm s(\bm u, \nabla \bm u)  \quad \mbox{in }\Omega .
\end{equation}
This general system of conservation laws encompasses a broad class of partial differential equations, including the heat equation, convection–diffusion problems, linear and nonlinear elasticity, and the compressible Euler and Navier–Stokes equations.

We denote by $\mathcal{T}_h$ a collection of disjoint regular elements $K$ that partition $\Omega$, and set $\partial \mathcal{T}_h := \{ \partial K : K \in \mathcal{T}_h \} $ to be the collection of the boundaries of the elements in $\mathcal{T}_h$. Let $\mathcal{F}_h$ be a collection of faces in $\mathcal{T}_h$.  Let $\mathcal{P}^{k}(K)$ denote the space of complete polynomials of degree $k$ on a domain $K \in \mathbb{R}^n$, let $L^2(K)$ be the space of square-integrable functions on $K$. We  introduce the following discontinuous finite element spaces:
$$\bm{\mathcal{Q}}_{h}^k  = \big\{\bm{q} \in [L^2(\mathcal{T}_h)]^{m \times d} \ : \ \bm q|_K \in \bm W(K),   \ \ \forall K \in \mathcal{T}_h \big\} , $$
$$\bm{\mathcal{V}}_{h}^k  = \big\{\bm{v} \in [L^2(\mathcal{T}_h)]^m \ : \ \bm v|_K \in \bm V(K), \ \ \forall K \in \mathcal{T}_h \big\} , $$
$$\bm{\mathcal{M}}_{h}^k  = \big\{\bm{\mu} \in [L^2(\mathcal{F}_h)]^m \ : \ \bm \mu|_F \in \bm V(F), \ \ \forall F \in \mathcal{F}_h \big\} , $$
where $\bm{W}(K) \equiv [\mathcal{P}^k(K)]^{m \times d}$ and $\bm{V}(K) \equiv [\mathcal{P}^k(K)]^{m}$. Next, we define several inner products associated with these finite element spaces as
\begin{subequations}
\label{innerProducts}
\begin{alignat}{3}
& (\bm{w},\bm{v})_{\mathcal{T}_h} && = \sum_{K \in \mathcal{T}_h} (\bm{w}, \bm{v})_K && = \sum_{K \in \mathcal{T}_h} \int_{K} \bm{w} \cdot \bm{v}  , \nonumber \\
& (\bm{q},\bm{p})_{\mathcal{T}_h} && = \sum_{K \in \mathcal{T}_h} (\bm{q}, \bm{p})_K && = \sum_{K \in \mathcal{T}_h} \int_{K} \bm{q} : \bm{p} , \nonumber \\
& \left\langle \bm{w}, \bm{v} \right\rangle_{\partial \mathcal{T}_h} && = \sum_{K \in \mathcal{T}_h} \left\langle \bm{w},\bm{v} \right\rangle_{\partial K} && = \sum_{K \in \mathcal{T}_h} \int_{\partial K} \bm{w} \cdot \bm{v} , \nonumber
\end{alignat}
\end{subequations}
for $\bm{w}, \bm{v} \in \bm{\mathcal{V}}_{h}^k$, $\bm{q}, \bm{p} \in \bm{\mathcal{Q}}_{h}^k$,  where $\cdot$ and $:$ denote the scalar product and Frobenius inner product, respectively. 

The HDG discretization of the governing equations reads as follows: Find $\big( \bm{q}_h,\bm{u}_h,  \widehat{\bm{u}}_h\big) \in \bm{\mathcal{Q}}_h^k \times \bm{\mathcal{V}}_h^k \times \bm{\mathcal{M}}_h^k$ such that
\begin{subequations}
\label{HDG}
\begin{alignat}{1}
\label{HDGa}
\big( \bm{q}_h, \bm{r} \big) _{\mathcal{T}_h} + \big( \bm{u}_h, \nabla \cdot \bm{r} \big)  _{\mathcal{T}_h} -  \big< \widehat{\bm{u}}_h, \bm{r} \cdot \bm{n} \big> _{\partial \mathcal{T}_h}  & =  0, \\
\label{HDGb}
\Big( \frac{\partial \, \bm{u}_h}{\partial t}, \bm{w} \Big)_{\mathcal{T}_h} - \Big(  \bm{F}(\bm u_h, \bm q_h), \nabla \bm{w} \Big) _{\mathcal{T}_h}  +  \left\langle \widehat{\bm{f}}_h, \bm{w} \right\rangle_{\partial \mathcal{T}_h}  & = \Big(  \bm{s}(\bm u_h, \bm q_h), \bm{w} \Big) _{\mathcal{T}_h} , \\
\label{HDGd}
  \left\langle \widehat{\bm{f}}_h, \bm{\mu} \right\rangle_{\partial \mathcal{T}_h\backslash \partial \Omega} + \left\langle \widehat{\bm{b}}_h, \bm{\mu} \right\rangle_{\partial \Omega}  & =  0, 
\end{alignat}
for all $\big( \bm{r},\bm{w}, \bm \mu\big) \in \bm{\mathcal{Q}}_h^k \times \bm{\mathcal{V}}_h^k \times \bm{\mathcal{M}}_h^k$, where the numerical flux $\widehat{\bm{f}}_h$ is defined as
\begin{equation}
\label{numflux2}
\widehat{\bm{f}}_h  = \bm{F}(\widehat{\bm u}_h , {\bm{q}}_h)    \cdot \bm{n} + \bm S(\bm u_h, \widehat{\bm u}_h) \cdot (\bm u_h - \widehat{\bm u}_h)
\end{equation}
\end{subequations}
and $\bm S(\bm u_h, \widehat{\bm u}_h)$ is the stabilization matrix. The boundary flux $\widehat{\bm{b}}_h$ is defined according to the boundary conditions. The HDG discretization of the steady-state system (\ref{eqst}) can be directly obtained from (\ref{HDG}) by omitting the time-derivative term.

The stabilization matrix can play an important role in the stability and accuracy of the HDG discretization. It governs the numerical fluxes exchanged across element boundaries and serves as a mechanism for introducing appropriate dissipation. Various choices for the stabilization matrix have been proposed, many of which are inspired by approximate Riemann solvers. Common schemes include the Roe solver~\cite{Roe1997}, the Lax–Friedrichs  flux~\cite{cockburn01:_rkdg}, and the HLL and HLLEM schemes~\cite{Harten83,Harten1997a,Vila-Perez2021}, with different levels of numerical dissipation and complexity. In the numerical experiments, we adopt the global Lax–Friedrichs stabilization due to its simplicity, ease of implementation, and robustness across a wide range of test cases.

Furthermore, the precise definition of the boundary flux  $\widehat{\bm{b}}_h$ in the HDG discretization (\ref{HDG}) depends on the boundary conditions. For Dirichlet boundary conditions, where the solution is prescribed as $\bm{u} = \bm{u}_D$ on $\partial \Omega$, the numerical boundary flux is defined as $\widehat{\bm{b}}_h = \widehat{\bm{u}}_h - \bm u_D$. For Neumann boundary conditions, where the normal component of the flux is specified as $\bm{q} \cdot \bm{n} = \bm{q}_N$ on $\partial \Omega$, the numerical boundary flux is defined by $\widehat{\bm{b}}_h = {\bm{q}}_h \cdot \bm n + \bm S(\bm u_h, \widehat{\bm{u}}_h) \cdot (\bm u_h -  \widehat{\bm{u}}_h) - \bm q_N$. More generally, for complex boundary conditions—such as mixed or Robin-type conditions, or in systems involving compressible flow with characteristic-based boundary treatments—the numerical boundary flux $\widehat{\bm{b}}_h$ must be defined accordingly to weakly enforce the correct physical behavior.

\subsection{Implicit Time Integration}

Next, we discretize the semi-discrete HDG formulation \eqref{HDG} in time using Diagonally Implicit Runge–Kutta (DIRK) schemes \cite{alexa77,Hairer1996}, which offer a flexible family of implicit time integration methods well-suited for stiff systems. A DIRK scheme with $s$ stages and formal order of accuracy $p$ is denoted as DIRK($s$, $p$). These schemes feature a lower triangular Butcher tableau with identical diagonal entries, allowing each stage to be solved implicitly but sequentially, which facilitates reuse of matrix structures and preconditioners between stages. The stage-by-stage structure of DIRK methods aligns naturally with the modularity of HDG solvers and enables efficient implementation on GPU architectures through reuse of assembled matrices, sparsity patterns, and Krylov subspace preconditioners. 

In this work, DIRK schemes are used for both transient and pseudo-transient  simulations. As a special case, the Backward-Euler method is a DIRK(1,1) scheme with a single stage and first-order accuracy. While it is only first-order accurate in time, it is unconditionally stable and often used as a robust time-stepping scheme in steady-state solvers, pseudo-transient continuation, or nonlinear homotopy methods. For unsteady problems requiring higher temporal accuracy, multi-stage DIRK schemes such as DIRK(2,2) or DIRK(3,3) are employed. These methods preserve the favorable stability properties of implicit integration while offering improved accuracy and flexibility in time-dependent simulations. 

Each DIRK($s$, $p$) scheme reduces to solving a sequence of $s$   problems that resemble DIRK(1,1)-like substeps, each involving a fully-discrete HDG system with updated right-hand sides. We describe the time integration procedure for the Backward-Euler scheme. Let $(\bm{q}_h^n,\bm{u}_h^n,\widehat{\bm{u}}_h^n,\widehat{\bm{f}}_h^n)$ denote the numerical approximations to $(\bm{q}_h(t^n),\bm{u}_h(t^n),\widehat{\bm{u}}_h(t^n),\widehat{\bm{f}}_h(t^n))$ at time level $t^n = n \Delta t$, where $\Delta t$ is the time step size. The Backward-Euler scheme advances the solution from $t^{n-1}$ to $t^{n}$ by solving the following nonlinear system: Find $\big( \bm{q}_h^n,\bm{u}_h^n,  \widehat{\bm{u}}_h^n\big) \in \bm{\mathcal{Q}}_h^k \times \bm{\mathcal{V}}_h^k \times \bm{\mathcal{M}}_h^k$ such that
\begin{subequations}
\label{HDGBE}
\begin{alignat}{1}
\label{HDGBEa}
\big( \bm{q}_h^n, \bm{r} \big) _{\mathcal{T}_h} + \big( \bm{u}_h^n, \nabla \cdot \bm{r} \big)  _{\mathcal{T}_h} -  \big< \widehat{\bm{u}}_h^n, \bm{r} \cdot \bm{n} \big> _{\partial \mathcal{T}_h}  & =  0, \\
\label{HDGBEb}
\Big( \frac{ \bm{u}_h^n}{\Delta t}, \bm{w} \Big)_{\mathcal{T}_h} - \Big(  \bm{F}(\bm u_h^n, \bm q_h^n), \nabla \bm{w} \Big) _{\mathcal{T}_h}  +  \left\langle \widehat{\bm{f}}_h^n, \bm{w} \right\rangle_{\partial \mathcal{T}_h}  & = \Big( \bm s(\bm u_h^n, \bm q_h^n) + \frac{  \bm{u}_h^{n-1}}{\Delta t}, \bm{w} \Big)_{\mathcal{T}_h}, \\
\label{HDGBEd}
  \left\langle \widehat{\bm{f}}_h^n, \bm{\mu} \right\rangle_{\partial \mathcal{T}_h\backslash \partial \Omega} + \left\langle \widehat{\bm{b}}_h^n, \bm{\mu} \right\rangle_{\partial \Omega}  & =  0, 
\end{alignat}
for all $\big( \bm{r},\bm{w}, \bm \mu\big) \in \bm{\mathcal{Q}}_h^k \times \bm{\mathcal{V}}_h^k \times \bm{\mathcal{M}}_h^k$, where the numerical flux $\widehat{\bm{f}}_h^n$ is given by
\begin{equation}
\label{numfluxBE2}
\widehat{\bm{f}}_h^n  = \bm{F}(\widehat{\bm u}_h^n , {\bm{q}}_h^n)    \cdot \bm{n} + \bm S(\bm u_h^n, \widehat{\bm u}_h^n) \cdot (\bm u_h^n - \widehat{\bm u}_h^n) .
\end{equation}
\end{subequations}
 Higher-order DIRK schemes follow a similar pattern with repeated application of the single-stage formulation across multiple stages. We refer to \cite{Nguyen2012} for a detailed description of the implementation of higher-order DIRK schemes with HDG discretization. 

\subsection{Linearization and Static Condensation}

For nonlinear PDEs, the fully discrete HDG system \eqref{HDGBE} is nonlinear and requires iterative solution techniques. By applying the Newton-Raphson method to linearize (\ref{HDGBE}) around the current solution $(\bm{q}_h^{ni}, \bm u_h^{ni}, \widehat{\bm u}_h^{ni})$ at iteration $i$, we obtain a linear system of the form
\begin{equation}
\label{matrixform} \left[
\begin{array}{ccc}
\mathbf{A}^{ni} & \quad \mathbf{B}^{ni} \quad  & \mathbf{C}^{ni} \\
\mathbf{D}^{ni} & \quad \mathbf{E}^{ni} \quad & \mathbf{F}^{ni}\\
\mathbf{G}^{ni} & \quad \mathbf{H}^{ni} \quad & \mathbf{J}^{ni}
\end{array}
\right] \left[
\begin{array}{c}
\delta \mathbf{q}^{ni} \\
\delta \mathbf{u}^{ni} \\
\delta \widehat{\mathbf{u}}^{ni}
\end{array}
\right] =\left[
\begin{array}{c}
\mathbf{r}^{ni}_{\bm q} \\
\mathbf{r}^{ni}_{\bm u} \\
\mathbf{r}^{ni}_{\widehat{\bm u}}
\end{array}
\right] .
\end{equation}
Here $(\delta \mathbf{q}^{ni}, \delta \mathbf{u}^{ni}, \delta \widehat{\mathbf{u}}^{ni})$
represents the vector of degrees of freedom for the solution increment $(\delta \bm{q}_h^{ni}, \delta \bm u_h^{ni}, \delta \widehat{\bm u}_h^{ni})$.  The next Newton iterate is computed via
\begin{equation}
(\bm{q}_h^{ni+1},  \bm u_h^{ni+1}, \widehat{\bm u}_h^{ni+1}) =  (\bm{q}_h^{ni},  \bm u_h^{ni}, \widehat{\bm u}_h^{ni}) + \alpha (\delta \bm{q}_h^{ni}, \delta \bm u_h^{ni}, \delta \widehat{\bm u}_h^{ni})
\end{equation}
where the constant $\alpha$ is chosen through line search to ensure monotonic decrease of the residual norm. The Newton-Raphson procedure terminates when the residual norm is less than a prescribed tolerance. 

Next, static condensation is employed to efficiently solve the linear system (\ref{matrixform}). It is done by eliminating $(\delta \mathbf{q}^{ni}, \delta \mathbf{u}^{ni})$ to obtain a reduced globally coupled system
 only for $\delta \widehat{\mathbf{u}}^{ni}$ as
\begin{subequations}
\label{linearsystem}
\begin{alignat}{1}
\label{eq13-a} &\mathbf{K}^{ni} \; \delta \widehat{\mathbf{u}}^{ni} = \mathbf{r}^{ni},
\intertext{where} \label{eq13-b} &\mathbf{K}^{ni}= {\mathbf{J}}^{ni} -\left[
\begin{array}{cc}
\mathbf{G}^{ni} & \mathbf{H}^{ni}
\end{array}
\right]\left[
\begin{array}{ccc}
\mathbf{A}^{ni} & \mathbf{B}^{ni} \\
\mathbf{D}^{ni} & \mathbf{E}^{ni} \\
\end{array}
\right]^{-1} \left[
\begin{array}{c}
\mathbf{C}^{ni} \\
\mathbf{F}^{ni} \\
\end{array}
\right], \intertext{and} \label{eq13-c} &\mathbf{r}^{ni} =
\mathbf{r}^{ni}_{\widehat{\bm u}} -\left[
\begin{array}{cc}
\mathbf{G}^{ni} & \mathbf{H}^{ni}
\end{array}
\right]\left[
\begin{array}{ccc}
\mathbf{A}^{ni} & \mathbf{B}^{ni} \\
\mathbf{D}^{ni} & \mathbf{E}^{ni} \\
\end{array}
\right]^{-1} \left[
\begin{array}{c}
\mathbf{r}^{ni}_{\bm q} \\
\mathbf{r}^{ni}_{\bm u}
\end{array}
\right].
\end{alignat}
\end{subequations}
Once the global linear system (\ref{eq13-a}) is solved, we compute $(\delta \mathbf{q}^{ni}, \delta \mathbf{u}^{ni})$ as
\begin{equation}
\left[
\begin{array}{c}
\delta \mathbf{q}^{ni} \\
\delta \mathbf{u}^{ni}
\end{array}
\right] = \left[
\begin{array}{ccc}
\mathbf{A}^{ni} & \mathbf{B}^{ni} \\
\mathbf{D}^{ni} & \mathbf{E}^{ni} \\
\end{array}
\right]^{-1} \left(
\left[
\begin{array}{c}
\mathbf{r}^{ni}_{\bm q} \\
\mathbf{r}^{ni}_{\bm u}
\end{array}
\right] -  \left[
\begin{array}{c}
\mathbf{C}^{ni} \\
\mathbf{F}^{ni} \\
\end{array}
\right] \delta \widehat{\mathbf{u}}^{ni}   
\right) .
\end{equation}
Due to the discontinuous nature of the approximation spaces, the matrix
$[\mathbf{A}^{ni} \ \mathbf{B}^{ni};  \mathbf{D}^{ni} \  \mathbf{E}^{ni}]$ is
block-diagonal and inverted in an element-by-element fashion to
yield a block-diagonal inverse. Therefore, the stiffness matrix $\mathbf{K}^{ni}$ and the residual vector $\mathbf{r}^{ni}$ can be assembled efficiently in parallel on GPUs.

\subsection{Structure of the Global Linear System}

The global linear system \eqref{eq13-a}, arising from the static condensation of the HDG formulation, exhibits a special sparse structure. The matrix $\mathbf{K}^{ni}$ couples only the degrees of freedom associated with the numerical trace $\widehat{\bm{u}}_h$ on the mesh skeleton, namely, the union of all element faces. This results in a global system whose size is significantly smaller than that of the original system involving both element and face degrees of freedom. Specifically, $\mathbf{K}^{ni}$ has dimensions $N_{\rm dof} \times N_{\rm dof}$, where  $N_{\rm dof} = M \times P_{F} \times N_F$, with $M$ denoting the number of state variables, $P_F$ the number of polynomial basis functions per face, and $N_F$  the total number of faces in the mesh.

Each entry of $\mathbf{K}^{ni}$ corresponds to the interaction between the trace unknowns on two neighboring faces.  Consequently, the matrix exhibits a sparse block structure governed by face connectivity: each face is coupled only to its neighboring faces through the shared elements. Specifically, only trace degrees of freedom on faces that belong to the same or adjacent elements are coupled. To represent this structure efficiently, $\mathbf{K}^{ni}$ is stored as a dense block format with dimensions $M \times P_F \times M \times P_F \times (2N_{\rm lfe} - 1) \times N_F$, where $N_{\rm lfe}$ is the number of local faces per element. The factor $(2N_{{\rm lfe}} - 1)$ arises because each face interacts with all other faces in the two elements sharing that face except for double-counting the central face. This localized coupling structure enables the use of block-wise memory layouts, block-based preconditioners, domain decomposition strategies, and GPU-based solvers.

\section{Iterative Linear Solvers and Preconditioning Techniques}

\subsection{Restarted GMRES}

To solve the global linear system \eqref{eq13-a}, we employ the Generalized Minimal Residual (GMRES) method \cite{sasc86}, a Krylov subspace iterative solver well suited for large, sparse, and nonsymmetric systems that arise from HDG discretizations. The convergence of GMRES can be significantly improved through preconditioning. Specifically, we apply left preconditioning and solve the equivalent system
\begin{equation}
\left(\mathbf{P}^{ni} \right)^{-1}  \mathbf{K}^{ni} \delta \widehat{\mathbf{u}}^{ni} = \left(\mathbf{P}^{ni} \right)^{-1}  \mathbf{r}^{ni},
\end{equation}
where $\left(\mathbf{P}^{ni} \right)^{-1}$ is a preconditioner matrix designed to approximate the inverse of $\mathbf{K}^{ni}$ in a computationally efficient manner. The construction of $\mathbf{P}^{ni}$ plays a central role in the efficiency and robustness of the iterative solver.

As GMRES must store and orthogonalize an increasingly large Krylov subspace,  the memory storage and orthogonalization cost increase with the number of iterations. We use the restarted GMRES denoted by GMRES($m$) to limit the dimension of the Krylov subspace to $m$ vectors. This reduces memory footprint and computational overhead, making it more suitable for GPU implementations where memory bandwidth and storage are limited. Typically, $m$ is chosen in the range of 20 to 100 to balance convergence rate, memory usage, and GPU parallelization efficiency.

\subsection{Block-Jacobi Preconditioner}

A highly parallelizable preconditioning strategy is the Block-Jacobi preconditioner. Due to the discontinuous nature of HDG discretizations, the global matrix $\mathbf{K}^{ni}$ exhibits a sparse, block-structured pattern in which only trace degrees of freedom on neighboring faces are coupled. The Block-Jacobi preconditioner $\mathbf{P}^{ni}$ is constructed by extracting the diagonal blocks of $\mathbf{K}^{ni}$ corresponding to each face. Hence, the inverse of $\mathbf{P}^{ni}$ is also a block-diagonal matrix which is stored as a dense block format with dimensions $M \times P_F \times M \times P_F \times N_F$. This preconditioner offers several key advantages. First, since each block can be inverted independently, the Block-Jacobi preconditioner is embarrassingly parallel and maps naturally to GPU architectures. Second, only the diagonal blocks need to be stored, reducing memory requirements compared to block ILU preconditioners. And third, the preconditioner requires no inter-face communication, making it particularly suitable for high-throughput computing. In our implementation, the block matrices  are extracted and stored in dense format, allowing for fast batched inversion using GPU-accelerated linear algebra libraries such as cuBLAS or hipBLAS. 

Block-Jacobi preconditioner provides a practical balance between simplicity, speed, and scalability for a wide range of problems. However, its effectiveness is limited by the absence of global coupling information across blocks, which can result in slow convergence for problems with strong cross-element interactions or poor conditioning. It is often used as a baseline preconditioner and serves as a building block for more sophisticated strategies.

\subsection{Additive Schwarz Preconditioner}

To improve upon the convergence limitation of Block-Jacobi while retaining parallelism, we consider  Additive Schwarz Method (ASM) \cite{Fernandez2017a}. The ASM preconditioner augments the purely local structure of Block-Jacobi by incorporating limited global coupling through overlapping subdomains. This overlap allows for better communication of information across the mesh, leading to significantly improved convergence rates especially for problems with strong inter-element coupling and ill-conditioned systems. Although ASM introduces more memory and computational overhead compared to Block-Jacobi, the improved convergence behavior often leads to lower total time to solution for large-scale problems.

We construct the ASM preconditioner by decomposing the mesh into a set of overlapping subdomains, where each subdomain corresponds to either a single element or a small patch of elements. Let $\mathbf{R}_{\ell}$ be a restriction operator that extracts the degrees of freedom associated with subdomain ${\ell}$, and let $\mathbf{K}_{\ell} = \mathbf{R}_{\ell} \mathbf{K}^{ni} \mathbf{R}_{\ell}^T$ denote the local matrix for subdomain ${\ell}$. The ASM preconditioner is then defined as
\begin{equation}
\left(\mathbf{P}^{ni} \right)^{-1} = \sum_{{\ell}=1}^{N_s} \mathbf{R}_{\ell}^T \mathbf{K}_{\ell}^{-1} \mathbf{R}_{\ell},
\end{equation}
where $N_s$ is the total number of overlapping subdomains. Several key aspects of this preconditioner make it well suited for GPU implementations. Each subdomain problem can be solved independently, enabling high degrees of parallelism. The overlap between subdomains enables the preconditioner to better approximate the global coupling structure, significantly accelerating convergence relative to Block-Jacobi. The subdomain size and degree of overlap can be tuned to balance memory usage, convergence speed, and computational workload. The local matrices $\mathbf{K}_{\ell}$ are small dense systems stored in a batched format and inverted using high-performance GPU libraries like cuBLAS and hipBLAS.

\subsection{Polynomial Precondioner}

As an alternative to factorization-based preconditioners, we consider polynomial preconditioners \cite{Loe2022}, which offer a matrix-free, highly parallelizable, and architecture-friendly approach to accelerating Krylov subspace solvers. Unlike ILU or Schwarz-type preconditioners, polynomial preconditioners avoid explicit factorizations or local solves by approximating the action of the matrix inverse through repeated matrix-vector products. This makes them especially attractive for GPU-based solvers, where parallel throughput and memory bandwidth are the primary performance constraints.

The polynomial preconditioner approximates the inverse of the global matrix $\mathbf{K}^{ni}$ by a truncated polynomial expansion [refs] as follows
\begin{equation}
\left(\mathbf{P}^{ni} \right)^{-1} = \sum_{j=1}^{P} \mathbf{Q}^{ni}_j, \quad \mbox{where }  \mathbf{Q}^{ni}_j = \left(\mathbf{I} - \frac{\mathbf{K}^{ni}}{\theta_1}  \right) \ldots \left(\mathbf{I} - \frac{\mathbf{K}^{ni}}{\theta_{j-1}}  \right)  \frac{1}{\theta_j} 
\end{equation}
where ${\theta_j}$ are harmonic Ritz values obtained by performing one cycle of GMRES($P$) with a random initial vector and extracting eigenvalue approximations from the associated Hessenberg matrix. The degree $P$ controls the quality of the approximation and introduces a trade-off between preconditioner cost and effectiveness. The application of the polynomial preconditioner to a vector involves repeated matrix-vector products with $\mathbf{K}^{ni}$, which are highly parallelizable and well suited for GPU acceleration. We implement this preconditioner using batched matrix-vector products and recurrence evaluations. The recurrence relation is fused with GMRES vector operations to maximize data reuse and minimize synchronization overhead.

This approach is particularly beneficial for ill-conditioned systems, where GMRES requires a large number of iterations, leading to expensive orthogonalization and memory usage. Polynomial preconditioner can substantially reduce the iteration count by accelerating convergence in such cases. Furthermore, polynomial preconditioners can be applied not only to the original system but also to a preconditioned system, allowing them to enhance other preconditioners. This hybrid strategy is especially effective in large-scale, high-order problems where no single preconditioner performs well across all spectral modes.

\section{GPU Implementation}

This section describes the GPU implementation of the HDG method with emphasis on the efficient realization of the local solvers, the assembly of the global linear system, and the Newton-GMRES solver with various preconditioning strategies. First, we describe the element-local kernels for eliminating interior degrees of freedom and computing local contributions in parallel across the mesh. Next, we outline the assembly of the condensed skeleton system in a block-sparse format suitable for GPU evaluation. We then present the GPU realization of the Newton-Krylov solver, focusing on the matrix–vector product and the integration of preconditioners such as Block-Jacobi, Additive Schwarz, and polynomial preconditioners.  Our goal is to exploit the high throughput and memory bandwidth of modern GPUs by reformulating all core kernels in terms of dense-block operations that map naturally to batched BLAS routines such as those provided by \texttt{cuBLAS} and \texttt{hipBLAS}. 

In the remainder of this paper we omit the superscript $ni$ when denoting matrices, vectors, and operators  at time step $n$ and Newton iteration $i$. Unless otherwise specified, all symbols should be understood as referring to the current iteration state. When distinction between different time levels or Newton iterations is essential, the corresponding superscripts will be explicitly reinstated.

\subsection{Implementing the Local Solvers}

It follows from the HDG weak formulation \eqref{HDGBE} that the degrees of freedom of $\bm q_h$ can be expressed locally in terms of those of $\bm u_h$ and $\widehat{\bm u}h$ as
\begin{equation}
\label{eqw15}
\mathbf{q}^e = - \left( \mathbf{A}^e \right)^{-1} \big( \mathbf{B}^e \mathbf{u}^e + \mathbf{C}^e \widehat{\mathbf{u}}^e \big),
\qquad \forall K^e \in \mathcal{T}h ,
\end{equation}
where $K^e$ denotes an element of the mesh $\mathcal{T}h$, and $\mathbf{A}^e$, $\mathbf{B}^e$, and $\mathbf{C}^e$ are elemental matrices arising from the discretization.
In expanded form, this relation can be written as
\begin{equation}
\label{eqw16}
\mathbf{q}^e_{dm} = - \left( \mathbf{M}^e \right)^{-1}  (\mathbf{B}^{e}_d \mathbf{u}^e_m + \mathbf{C}^{e}_d \widehat{\mathbf{u}}^e_m ), \quad 1 \le m \le M, 1 \le d \le D, 
\end{equation}
where the matrices $\mathbf{M}^e$, $\mathbf{B}^e_d$, and $\mathbf{C}^e_d$ are defined by
\begin{equation}
\begin{split}
M^e_{ij} & = \int_{K^e} \phi_i^e(\bm x) \phi^e_j(\bm x)  d \bm x, \\
B^{e}_{dij} & = \int_{K^e} \phi_i^e(\bm x) \frac{\partial \phi^e_j(\bm x)}{\partial x_d}  d \bm x, \\
C^{e}_{dil} & = -\int_{\partial K^e} \phi_i^e(\bm x) \psi^e_l(\bm x) n_d^e  d \bm x ,
\end{split}
\end{equation}
for $1 \le i,j \le P_E$ and $1 \le l \le N_{{\rm lfe}} P_F$, with $P_E$ and $P_F$ denoting the number of basis functions on each element and face, respectively, and $N_{{\rm lfe}}$ the number of local faces per element. Here $\phi^e_i$ are basis functions of $\mathcal{P}^k(K^e)$ and $\psi^e_l$ are basis functions of $\mathcal{P}^k(\partial K^e)$. 

The element mass matrices $\mathbf{M}^e$ for all elements $1 \le e \le N_E$ can be computed simultaneously in matrix form as
\begin{equation}
\mathbf{M} = [\mathbf{M}^1, \ldots, \mathbf{M}^{N_E}] = \bm{\Phi} \bm{W},
\end{equation}
where $\bm{\Phi} \in \mathbb{R}^{(P_E^2) \times Q_E}$ stores the values of the products of the basis functions  evaluated at the quadrature points on the reference element, and $\bm{W} \in \mathbb{R}^{Q_E \times N_E}$ collects the corresponding quadrature weights multiplied by the Jacobian determinants of all elements. Here $Q_E$ denotes the number of quadrature points used on the reference element. The matrices $\mathbf{B}^e_d$ and $\mathbf{C}^e_d$ are computed in a similar manner.  The inverses $(\mathbf{M}^e)^{-1}$ are obtained using batched dense matrix inversion routines provided by \texttt{cuBLAS} or \texttt{hipBLAS}. The matrices $(\mathbf{M}^e)^{-1} \mathbf{B}^e_d$ and $(\mathbf{M}^e)^{-1} \mathbf{C}^e_d$ are precomputed once during setup and stored in GPU memory. During the solve phase, the evaluation of \eqref{eqw16} is performed via batched matrix–matrix multiplications, also using optimized BLAS kernels. This strategy allows the computation of $\bm q_h$ to be carried out efficiently across all elements on the GPU.

It follows from (\ref{eqw15}) that $\mathbf{r}^e_{\bm q} = 0$.  Hence, the first row of \eqref{matrixform} can be eliminated and substituted into the second  row to yield
\begin{equation}
\left(\mathbf{E}^e - \mathbf{D}^e  \left( \mathbf{A}^e \right)^{-1} \mathbf{B}^e \right) \delta \mathbf{u}^e + \left(\mathbf{F}^e - \mathbf{D}^e  \left( \mathbf{A}^e \right)^{-1} \mathbf{C}^e \right) \delta \widehat{\mathbf{u}}^e = \mathbf{r}^e_{\bm u}   
\end{equation}
which is equivalent to
\begin{equation}
\label{eqw19}
\delta \mathbf{u}^e = \left(\bar{\mathbf{E}}^e \right)^{-1} \left(\mathbf{r}^e_{\bm u} - \bar{\mathbf{F}}^e  \delta \widehat{\mathbf{u}}^e \right) 
\end{equation}
where $\bar{\mathbf{E}}^e = \mathbf{E}^e - \mathbf{D}^e  \left( \mathbf{A}^e \right)^{-1} \mathbf{B}^e$ and $\bar{\mathbf{F}}^e = \mathbf{F}^e - \mathbf{D}^e  \left( \mathbf{A}^e \right)^{-1} \mathbf{C}^e$. The matrices $\mathbf{D}^e, \mathbf{E}^e, \mathbf{F}^e$ and vectors $\mathbf{r}^e_{\bm u}$ are derived from the linearization of the equation (\ref{HDGBEb}). They are computed by evaluating fluxes, source terms, and their derivatives at quadrature points for all elements at once and multiplying the resulting arrays with the  basis function matrices via batched matrix–matrix multiplication BLAS kernels provided by \texttt{cuBLAS} or \texttt{hipBLAS}. Specifically, $\mathbf{E}^e \in \mathbb{R}^{P_E M \times P_E M }$ are given by
\begin{equation}
{E}^e_{iajb} = {E}^{K^e}_{iajb} + {E}^{\partial K^e}_{iajb}, \qquad 1 \le i, j \le P_E, 1 \le a, b \le M,    
\end{equation}
where 
\begin{equation}
\begin{split}    
{E}^{K^e}_{iajb}  & = \Big(\delta_{ab} \frac{ \phi_j}{\Delta t} - \frac{\partial s_a(\bm u_h, \bm q_h)}{\partial u_b} \phi_j, \phi_i \Big)_{K^e} - \Big(  \frac{\partial F_a(\bm u_h, \bm q_h)}{\partial u_b} \phi_j, \nabla \phi_i \Big) _{K^e}  \\
{E}^{\partial K^e}_{iajb}  & =  \left\langle \frac{\partial \widehat{f}_{h,a}(\bm u_h, \bm q_h, \widehat{\bm u}_h)}{\partial u_b} \phi_j, \phi_i \right\rangle_{\partial K^e}   .
\end{split}
\end{equation}
The elemental matrices $\mathbf{E}^{K^e}$ for all elements $1 \le e \le N_E$ can be obtained simultaneously by computing
\begin{equation}
\mathbf{Q} = \bm{\Psi} \mathbf{U},
\end{equation}
and permuting $\mathbf{Q} \in \mathbb{R}^{P_E^2 \times M^2  \times N_E}$, where $\bm{\Psi} \in \mathbb{R}^{P_E^2 \times Q_E (D+1)}$ stores the values of the products of the basis functions and their derivatives  evaluated at the quadrature points on the reference element, and $\bm{U} \in \mathbb{R}^{Q_E (D+1) \times M^2  N_E}$ collects the  quadrature weights multiplied by the Jacobian determinants and the partial derivatives of the source term and fluxes at quadrature points for all physical elements. Similarly, the elemental matrices $\mathbf{E}^{\partial K^e}, 1 \le e \le N_E,$ are obtained by computing integrals on the element boundaries. 

This local elimination step is a key ingredient of the static condensation procedure: by expressing $\delta \mathbf{u}^e$ in terms of $\delta \widehat{\mathbf{u}}^e$, the global coupled system can be reduced to a skeleton-only system involving only the trace unknowns. Once $\delta \widehat{\mathbf{u}}$ is solved globally, the element-level updates $\delta \mathbf{u}^e$ (and subsequently $\delta \mathbf{q}^e$) can be recovered locally in parallel across elements. This separation between local elimination and global trace solve is what enables HDG methods to combine reduced global system size with high parallel efficiency, particularly on GPU architectures.

\subsection{Assembling the Global Linear System}

Substituting the first row of \eqref{matrixform} into the third row yields
\begin{equation}
\left(\mathbf{H}^e - \mathbf{G}^e  \left( \mathbf{A}^e \right)^{-1} \mathbf{B}^e \right) \delta \mathbf{u}^e + \left(\mathbf{J}^e - \mathbf{G}^e  \left( \mathbf{A}^e \right)^{-1} \mathbf{C}^e \right) \delta \widehat{\mathbf{u}}^e = \mathbf{r}^e_{\widehat{\bm u}}   
\end{equation}
which can be written more compactly as
\begin{equation}
\label{eqw21}
\bar{\mathbf{H}}^e  \delta \mathbf{u}^e + \bar{\mathbf{J}}^e  \delta \widehat{\mathbf{u}}^e = \mathbf{r}^e_{\widehat{\bm u}}   
\end{equation}
where $\bar{\mathbf{H}}^e = \mathbf{H}^e - \mathbf{G}^e  \left( \mathbf{A}^e \right)^{-1} \mathbf{B}^e$ and $\bar{\mathbf{J}}^e = \mathbf{J}^e - \mathbf{G}^e  \left( \mathbf{A}^e \right)^{-1} \mathbf{C}^e$. 
The matrices $\mathbf{G}^e, \mathbf{H}^e, \mathbf{J}^e$ and the residual $\mathbf{r}^e_{\widehat{\bm u}}$ arise from the linearization of the numerical flux equation \eqref{HDGBEd}. They are assembled by evaluating the numerical and boundary fluxes, together with their derivatives, at quadrature points on all element faces. The resulting arrays are then multiplied by the face shape function matrices to obtain the final contributions. In practice, these operations are carried out in parallel on GPUs using batched matrix–matrix multiplication kernels from \texttt{cuBLAS} or \texttt{hipBLAS}.  

Next, substituting the local expression for $\delta \mathbf{u}^e$ from \eqref{eqw19} into \eqref{eqw21} yields the condensed elemental trace system
\begin{equation}
\label{eqw22}
\bar{\mathbf{K}}^e  \delta \widehat{\mathbf{u}}^e = \bar{\mathbf{r}}^e_{\widehat{\bm u}}, 
\end{equation}
where
\begin{equation}
\bar{\mathbf{K}}^e  = \bar{\mathbf{J}}^e  - \bar{\mathbf{H}}^e  \left(\bar{\mathbf{E}}^e \right)^{-1} \bar{\mathbf{F}}^e, \quad  \bar{\mathbf{r}}^e_{\widehat{\bm u}} = \mathbf{r}^e_{\widehat{\bm u}}  - \bar{\mathbf{H}}^e  \left(\bar{\mathbf{E}}^e \right)^{-1} \mathbf{r}^e_{\bm u} .
\end{equation}
The elemental stiffness matrices $\bar{\mathbf{K}}^e$ are not invertible on their own, since the trace degrees of freedom are globally coupled. Instead, they serve as the building blocks of the global condensed system \eqref{linearsystem}.  

The elemental matrices $\bar{\mathbf{K}}^e$ are stored in GPU memory as a dense-block array of dimensions $(M \times P_F \times N_{\mathrm{\rm lfe}})^2 \times N_E$. 
To assemble the global system, we traverse the mesh faces. For each face, we identify the two adjacent elements, extract the corresponding face blocks from their local matrices $\bar{\mathbf{K}}^e$, and merge them to form the global block contribution of size $(M \times P_F)^2 \times (2N_{\mathrm{\rm lfe}}-1)$. These face-level blocks are then inserted into the global matrix $\mathbf{K}$, which is stored in a block compressed row format of dimensions $(M \times P_F)^2 \times (2 N_{\mathrm{\rm lfe}} - 1) \times N_F$. The global residual vector $\mathbf{r}$ is formed in a similar manner and stored as a dense-block array of dimensions $M \times P_F \times N_F$.  This face-based assembly strategy exploits the natural sparsity pattern of the HDG skeleton system and is highly parallelizable. On GPUs, each face is processed independently by a warp or thread block, enabling simultaneous accumulation of contributions across all faces. By relying on dense-block representations and batched BLAS kernels, the implementation achieves both memory coalescence and high arithmetic intensity to ensure scalability. 

\subsection{Implementing Matrix-Vector Product}

Once the global condensed system \eqref{linearsystem} has been assembled, the dominant operation within Krylov solvers such as GMRES is the repeated evaluation of the matrix–vector product 
\[
\mathbf{y} = \mathbf{K} \, \mathbf{x},
\]
where $\mathbf{x}$ and $\mathbf{y}$ are global vectors of dimensions $M \times P_F \times N_F$. The vector $\mathbf{x}$ is first reshaped into an extended block format of size $M \times P_F (2N_{\mathrm{\rm lfe}} - 1) \times N_F$. This transformation ensures that, for each face, the degrees of freedom on the central face as well as those of its neighboring faces are stored contiguously in memory.   With this layout, the matrix–vector product reduces to a batched dense matrix–vector multiplication at the face level:
\[
\mathbf{y}_f = \mathbf{K}_f \, \mathbf{x}_f, \qquad f = 1, \ldots, N_F,
\]
where $\mathbf{K}_f \in \mathbb{R}^{(M P_F) \times (M P_F (2N_{\mathrm{\rm lfe}}-1))}$ denotes the block associated with face $f$, $\mathbf{x}_f$ is the corresponding local slice of the extended input vector, and $\mathbf{y}_f$ is the resulting contribution to the output vector.  The batched operations are carried out using \texttt{cuBlas} on NVIDIA GPUs or \texttt{hipBlas} on AMD GPUs. The strided-batched interface allows all $N_F$ face-level products to be executed in a single kernel launch, ensuring coalesced access, high parallelism and minimizing launch overhead. 

This implementation eliminates the need for sparse matrix storage formats and instead leverages the  block structure of HDG to exploit optimized dense linear algebra kernels. As a result, the matrix–vector product achieves both high arithmetic intensity and excellent scalability, making it well-suited for iterative solvers on modern GPU architectures.

\subsection{Implementing the Preconditioned Newton-GMRES Solver}

The nonlinear system arising at each time step is solved using a damped Newton iteration, where each Newton correction is computed with a preconditioned GMRES solver. The nonlinear solver is implemented in the routine \texttt{NonlinearSolver}, which begins by computing the initial solution and residual norms. At each Newton step, the residual vector $\mathbf{R}(u)$ is assembled via a GPU kernel \texttt{AssembleResidual}, which evaluates both element and face flux contributions in parallel. The norm of the residual is monitored to determine convergence. For robustness, a damped Newton strategy is employed: after computing the update $\delta \widehat{\mathbf{u}}$, the solution is tentatively updated as $\widehat{\mathbf{u}} \leftarrow \widehat{\mathbf{u}} + \alpha \delta \widehat{\mathbf{u}}$, where $\alpha$ is reduced by successive halving if the residual norm fails to decrease. This ensures global convergence even for strongly nonlinear problems.

At each Newton step, the correction $\delta \widehat{\mathbf{u}}$ is obtained by solving the linear system $\mathbf{K} \, \delta \widehat{\mathbf{u}} = {\mathbf{r}}$ using the preconditioned GMRES solver. The main kernels in GMRES correspond to (i) matrix–vector product already described in the previous subsection, (ii) preconditioner application, and (iii) orthogonalization.  The solver constructs an orthonormal Krylov basis using either the modified Gram–Schmidt (\texttt{MGS}) or classical Gram–Schmidt (\texttt{CGS}) process, both of which are implemented with GPU-accelerated BLAS dot products and axpy operations. To maintain numerical stability, the small dense Hessenberg system that arises in the Arnoldi process is updated on the CPU: Givens rotations are applied to introduce zeros below the diagonal, and the least-squares problem is solved by back substitution. Convergence is checked at each iteration based on the relative residual norm. Upon convergence, the update $\delta \widehat{\mathbf{u}}$ is obtained by multiplying the Krylov basis vectors with the computed least-squares coefficients via batched matrix–vector multiplication BLAS kernels on the GPU.  



The Krylov vectors are preconditioned prior to orthogonalization. To this end, we consider three classes of preconditioners that exploit the block structure of the HDG formulation and map efficiently to GPU architectures:

\begin{itemize}
    \item \textbf{Block-Jacobi Preconditioner:} A simple yet parallel baseline, constructed by inverting the diagonal face blocks of $\mathbf{K}$ independently. Its embarrassingly parallel nature makes it well-suited for GPUs.
    \item \textbf{Additive Schwarz Method (ASM) Preconditioner:} An overlapping domain decomposition strategy where each element-local system is enriched with contributions from its immediate neighbors. This captures more of the global coupling while maintaining high concurrency.
    \item \textbf{Polynomial Preconditioners:} A strictly matrix-free alternative that approximates $\mathbf{K}^{-1}$ using a truncated polynomial expansion based on harmonic Ritz values. This approach replaces factorizations with repeated matrix–vector products, which are naturally GPU-friendly.
\end{itemize}

These preconditioners and their hybrid variants are implemented in the Newton–GMRES solver. The following subsections describe each preconditioner in detail.

 
\subsection{Implementing the Block-Jacobi Preconditioner}

The Block-Jacobi preconditioner exploits the block structure of the HDG condensed matrix $\mathbf{K}$ by inverting only the diagonal face blocks. Each block corresponds to the interaction of the trace degrees of freedom on a single face with itself, and has dimensions $(M P_F) \times (M P_F)$. Since these blocks are independent of one another, their inverses can be computed and applied in parallel across all faces. We extract the diagonal blocks $\mathbf{K}_{f,f} \in \mathbb{R}^{(M P_F) \times (M P_F)}$ for all faces $f=1,\ldots,N_F$ from the global block structure. These blocks are stored as a dense array of dimensions $(M P_F)^2 \times N_F$ in GPU memory. Each block is inverted once using batched dense matrix factorization and inversion routines provided by \texttt{cuBLAS} or \texttt{hipBLAS}. The resulting inverses $\mathbf{P}_f^{-1}$ are stored for use throughout the GMRES iterations.

During GMRES iterations, the application of the preconditioner reduces to a batched dense matrix–vector multiplication
\[
\mathbf{z}_f = \mathbf{P}_f^{-1} \, \mathbf{y}_f, \qquad f = 1,\ldots,N_F,
\]
where $\mathbf{y}_f \in \mathbb{R}^{M P_F}$ is the local slice of the global vector $\mathbf{y}$ corresponding to face $f$, and $\mathbf{z}_f$ is the output. The operation is executed in parallel across all faces using strided-batched \texttt{gemv} kernels. Since each block is small and uniform in size, the GPU achieves high occupancy and sustained throughput.

The Block-Jacobi preconditioner is simple to implement and embarrassingly parallel, making it highly efficient on GPUs. It requires no inter-face communication and leverages optimized dense linear algebra kernels, resulting in very low overhead. Its main limitation is the absence of global coupling information, which may slow convergence for problems with strong inter-element interactions. Nevertheless, its scalability and simplicity make Block-Jacobi an effective baseline preconditioner.

\subsection{Implementing the ASM Preconditioner}

The ASM preconditioner extends the Block-Jacobi strategy by introducing limited overlap between subdomains. In the context of HDG, a natural choice of subdomain is the element, which couples together the face unknowns on its $N_{\mathrm{\rm lfe}}$ local faces. The ASM preconditioner is constructed by forming element-level matrices $\bar{\mathbf{P}}^e$ that approximate the contribution of each element and its immediate neighbors.  

The ASM matrices $\bar{\mathbf{P}}^e$ are initialized by copying the elemental stiffness matrices $\bar{\mathbf{K}}^e$ of dimensions $(M P_F N_{\mathrm{\rm lfe}}) \times (M P_F N_{\mathrm{\rm lfe}}) \times N_E$. Next, for each interior face, we identify the two neighboring elements $K^{e_1}$ and $K^{e_2}$ that share the face, as well as the corresponding local face indices $l_1$ on $K^{e_1}$ and $l_2$ on $K^{e_2}$. The block rows and block columns of $\bar{\mathbf{K}}^e$ are naturally partitioned into $N_{\mathrm{\rm lfe}} \times N_{\mathrm{\rm lfe}}$ blocks, each of size $(M P_F) \times (M P_F)$, corresponding to couplings between local faces. The ASM update then enriches the diagonal blocks associated with the shared face as follows:
\[
\bar{\mathbf{P}}^{e_1}_{l_1l_1} \; \leftarrow \; \bar{\mathbf{K}}^{e_1}_{l_1l_1} + \bar{\mathbf{K}}^{e_2}_{l_2l_2}, 
\qquad
\bar{\mathbf{P}}^{e_2}_{l_2l_2} \; \leftarrow \; \bar{\mathbf{K}}^{e_1}_{l_1l_1} + \bar{\mathbf{K}}^{e_2}_{l_2l_2} .
\]
As a result, each $\bar{\mathbf{P}}^e$ incorporates both intra-element couplings and information from its immediate neighbors, thereby capturing a richer approximation of the global system. On GPUs, this update is implemented as a face-parallel kernel: each thread processes one degree of freedom on an interior face, locates the two adjacent elements and their local face indices, and accumulates the neighbor contributions into the appropriate entry of $\bar{\mathbf{P}}^e$. This face-based update completely avoids race conditions. 
Once assembled, the inverses $(\bar{\mathbf{P}}^e)^{-1}$ are computed using batched LU routines from \texttt{cuBLAS} or \texttt{hipBLAS}, and stored as a dense array of dimensions $(M P_F N_{\mathrm{\rm lfe}}) \times (M P_F N_{\mathrm{\rm lfe}}) \times N_E$.

During GMRES iterations, the preconditioner action $\mathbf{z} = ASM(\mathbf{y})$ is carried out as follows. First, we traverse each element $K^e$ and its local faces to construct $\mathbf{y}^e$ from $\mathbf{y}$. For each local face $l$ of element $K^e$, we identify the corresponding global face index $f(e,l)$. The slice of $\mathbf{y}$ associated with this global face is then copied into the $l$-th block of $\mathbf{y}^e$. Repeating this process for all local faces and all elements results in
$\mathbf{y}^e$ as a dense array of size $M \times P_F \times N_{\mathrm{\rm lfe}} \times N_E$. Next, we perform a batched matrix–vector multiplication at the element level:
\[
\mathbf{z}^e = (\bar{\mathbf{P}}^e)^{-1} \mathbf{y}^e, \qquad e = 1,\ldots,N_E,
\]
which is processed in parallel using strided-batched matrix-vector product kernels. Finally, we reassemble the global vector $\mathbf{z} \in \mathbb{R}^{M \times P_F \times N_F}$ from $\mathbf{z}^e$ by reversing the restriction process. For each element $K^e$ and local face $l$, the local contribution $\mathbf{z}^e_l \in \mathbb{R}^{M \times P_F}$ is scattered back to the slice of $\mathbf{z}$ associated with the global face index $f(e,l)$. Since a global face is shared by two neighboring elements, contributions from both elements are accumulated. This restriction–solve–prolongation cycle is implemented entirely on the GPU. Restriction and prolongation are executed as face-parallel kernels to ensure coalesced memory access, while the local solves are handled by batched BLAS routines.

The ASM preconditioner strikes a balance between scalability and convergence. Compared to Block-Jacobi, it incorporates limited global coupling by overlapping neighboring elements, which accelerates convergence for ill-conditioned systems and convection-dominated problems. The use of dense-block storage and batched GPU kernels ensures high throughput and efficient use of device memory. 

\subsection{Implementing the Polynomial Preconditioner}

We implement the polynomial preconditioner \cite{Loe2022} in a matrix-free fashion, so that its application reduces to a short recurrence involving repeated matrix–vector products. Given a degree $P$ of the polynomial preconditioner, the first step is to compute a set of harmonic Ritz values $\{\theta_j\}_{j=1}^P$, which define the recurrence coefficients. To this end, we perform one cycle of Arnoldi/GMRES of dimension $P$ on a random initial vector. The Arnoldi process is carried out by applying repeated matrix–vector products with $\mathbf{K}$, followed by modified Gram–Schmidt (MGS) orthogonalization. This yields an upper Hessenberg matrix $\bm H_{P+1,P}$ that captures the spectral behavior of the global matrix $\mathbf{K}$. The harmonic Ritz values are then obtained as the eigenvalues of $\bm H_{P,P} + H_{P+1,P}^2 (\bm H^T_{P,P})^{-1} \bm e_P \bm e_P^T$ with elementary coordinate vector $\bm e_P = (0,\ldots, 1)$ [refs]. The resulting Ritz values may be real or appear in complex conjugate pairs. To ensure numerical stability during the polynomial evaluation, the values are reordered using a Leja sequence, which improves robustness for clustered eigenvalues. Conjugate symmetry is explicitly enforced, so that the final polynomial preconditioner has real coefficients.

Once the Ritz values $\{\theta_j\}$ have been computed and reordered in a Leja sequence, the polynomial preconditioner is applied to a vector using a recurrence that mimics Horner’s rule but explicitly accounts for real and complex conjugate pairs. Specifically, given an input vector $\mathbf{y}$, the preconditioner action is evaluated via a Horner-like product of shifted operators
\[
\mathbf{z} 
= \sum_{j=1}^{P} \mathbf{Q}_j \mathbf{y}
= \sum_{j=1}^{P} \left[ \left(\mathbf{I} - \frac{\mathbf{K}}{\theta_1}\right) \!\cdots\!
\left(\mathbf{I} - \frac{\mathbf{K}}{\theta_{j-1}}\right) \frac{1}{\theta_j} \right] \mathbf{y}.
\]
In practice we apply the equivalent nested form that requires only $P$ applications of $\mathbf{K}$ and $P$ axpy operations. Two device buffers are used in ping–pong fashion to avoid extra allocations. If the Ritz value $\theta_j$ is real, the recurrence applies a single shifted operator:
\[
\mathbf{w} \;\leftarrow\; \mathbf{w} + \frac{1}{\theta_j} \mathbf{q}, 
\qquad
\mathbf{q} \;\leftarrow\; \mathbf{q} - \frac{1}{\theta_j} \mathbf{K}\mathbf{q}.
\]
Here $\mathbf{q}$ is the current working vector, $\mathbf{K}\mathbf{q}$ is obtained from the batched HDG matrix–vector product, and $\mathbf{w}$ accumulates the final preconditioned result. When the Ritz values appear as a conjugate pair $\theta_j = a+ib$ and $\theta_{j+1} = a-ib$, the recurrence is applied in paired form to ensure real arithmetic:
\[
\mathbf{w} \;\leftarrow\; \mathbf{w} + \frac{1}{a^2+b^2}\left( 2a \mathbf{q} - \mathbf{K}\mathbf{q} \right),
\qquad
\mathbf{q} \;\leftarrow\; \mathbf{q} - \frac{1}{a^2+b^2} \mathbf{K}^2 \mathbf{q}.
\]
This update preserves stability while avoiding explicit complex computations. The axpy and scaling operations in the recurrence are fused into custom kernels that operate directly on the face-major layout to maximize cache reuse and minimize kernel-launch overhead. 

This implementation makes the polynomial preconditioner entirely matrix-free: it requires only matrix–vector products and simple vector updates, both of which are highly parallel and GPU-friendly. The explicit handling of real and complex conjugate pairs ensures numerical stability, while Leja ordering provides robustness in the presence of clustered Ritz values. The cost of applying a degree-$P$ polynomial is $P$ matrix-vector products plus $O(P)$ axpy/scale operations. Storage consists of the working buffers (two vectors) and the coefficients $\{\theta_j\}$.


The polynomial preconditioner may be applied as a left preconditioner to the original system or composed with an existing preconditioner $\mathbf{M}^{-1}$ (e.g., Block–Jacobi or ASM) to form a hybrid scheme:
\[
\mathbf{z} = \sum_{j=1}^P \mathbf{Q}_j \big(\mathbf{M}^{-1}\mathbf{y}\big) .
\]
 Because the coefficients $\{\theta_j\}$ are fixed during a GMRES cycle, the method is compatible with standard GMRES.

\section{Numerical Experiments}

In this section, we evaluate the performance, scalability, and robustness of four GPU-based preconditioning strategies:
\begin{enumerate}
\item \textbf{Block-Jacobi (BJ)}: A purely local preconditioner based on face-block diagonal inversion; 
\item \textbf{Additive Schwarz Method (ASM)}: A  preconditioner employs one-element overlapping subdomains; 
\item \textbf{Block-Jacobi with PP (BJ-PP($P$))}: A hybrid preconditioner that combines BJ with polynomial preconditioner with degree $P$; 
\item \textbf{ASM with PP (ASM-PP($P$))}: A hybrid preconditioner that combines ASM with polynomial preconditioning with degree $P$.
\end{enumerate}


The numerical experiments encompass a broad range of PDE models to assess the effectiveness and robustness of each preconditioner across diverse physical regimes, mesh configurations, and polynomial orders. The evaluation focuses on convergence behavior, computational throughput, and overall runtime performance on GPU architectures. All computations are performed in double precision, and GPU timings are measured using explicit device synchronization to ensure accurate wall-clock profiling. The HDG discretization and solver framework is implemented in the open-source platform \texttt{Exasim}~\cite{VilaPerez2022}. All source codes and scripts used in the present study are publicly available at \url{https://github.com/exapde/Exasim}.

The solver parameters used for all the test cases are summarized as follows. The stabilization matrix is chosen as $\boldsymbol{S} = \tau \boldsymbol{I}$ with $\tau$ being chosen to provide a suitable balance between numerical dissipation and stability within the HDG formulation. The resulting nonlinear system is solved using a Newton–GMRES algorithm, with a nonlinear convergence tolerance of $10^{-8}$ and a linear solver tolerance of $10^{-6}$. Each linearized system arising at a Newton iteration is solved using a restarted GMRES method with a restart dimension of 50 and classical Gram–Schmidt orthogonalization. The maximum number of GMRES iterations is set to 1000. 

\subsection{Two-Dimensional Viscous Burgers Problem}

\begin{figure}[htbp]
	\centering
 \includegraphics[width=0.7\textwidth]{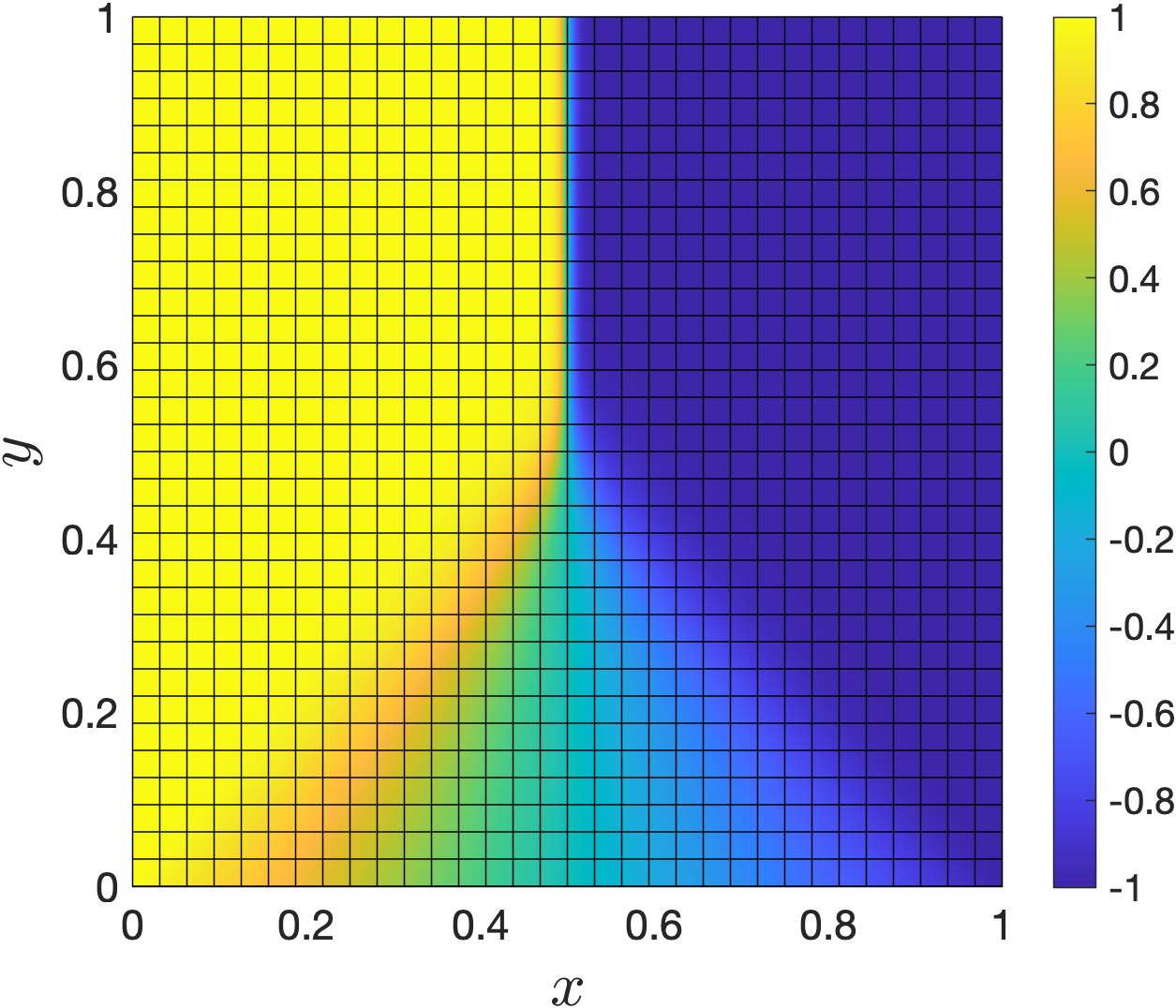}
 \caption{The numerical solution of the viscous Burgers problem for $k=3$ and $h = 1/32$.}
	\label{fig0}
\end{figure}

We consider the two-dimensional steady-state viscous Burgers problem \cite{Moro2012} in a unit square domain
\begin{equation}
\frac{\partial u}{\partial y}  +  \frac{1}{2} \frac{\partial u^2}{\partial x}  - \frac{1}{200} \left( \frac{\partial^2 u}{\partial x^2}  + \frac{\partial^2 u}{\partial y^2} \right) = 0, \quad \mbox{in } \Omega \equiv (0,1)^2    
\end{equation}
with Dirichlet boundary conditions $u(x,y) = 1 - 2x$ on $x = 0$, $x=1$, and $y=0$, while an outflow (Neumann-type) boundary condition is imposed on the top boundary $y=1$.  The computational domain is discretized using a sequence of structured quadrilateral meshes with uniform element sizes. The mesh resolution is characterized by the parameter $1/h \in 
\{16, 32, 64, 128\}$, corresponding to meshes with $16 \times 16$, $32 \times 32$, $64 \times 64$, and $128 \times 128$ elements, respectively. Each element supports a local polynomial approximation of degree $k \in \{1, 2, 3, 4\}$. Figure~\ref{fig0} shows the numerical solution computed using polynomial degree $k=3$ and mesh size $h = 1/32$. The solution remains smooth away from the centerline but exhibits a sharp transition near $x=0.5$, which is well resolved by the high-order HDG discretization.

Table~\ref{tab1} reports the computational time $t_{\mathrm{comp}}$ (in seconds) and the number of GMRES iterations $n_{\mathrm{gmres}}$ required to reach convergence for various polynomial orders $k$ and mesh resolutions $h$, comparing four preconditioning strategies. Since the degree of the polynomial preconditioner is fixed at 10, each GMRES iteration in BJ-PP(10) (respectively, ASM-PP(10)) requires ten successive applications of the BJ (respectively, ASM) preconditioner. The reported GMRES iteration counts correspond to the total number of linear iterations accumulated over all Newton iterations required to converge the steady-state problem. In other words, $n_{\mathrm{gmres}}$ represents the cumulative effort of the linear solver across the entire Newton process. Overall, the ASM-based preconditioners yield faster convergence and shorter computational times compared to the purely local BJ variants. In particular, ASM-PP(10) achieves the lowest iteration counts across all discretization levels, demonstrating the effectiveness of combining domain decomposition and polynomial acceleration. The computational cost increases with both polynomial degree $k$ and mesh refinement $1/h$, as expected, but the preconditioned solvers exhibit robust performance even on the finest grid considered ($1/h = 128$).

\begin{table}[htbp]
\centering
\begin{tabular}{|c|c|cc|cc|cc|cc|}
\hline
  \multicolumn{2}{|c|}{} & \multicolumn{2}{c|}{BJ} & \multicolumn{2}{c|}{BJ-PP(10)} & \multicolumn{2}{c|}{ASM} & \multicolumn{2}{c|}{ASM-PP(10)}   \\
\hline
$k$ & $1/h$ & $t_{\mathrm{comp}}$
  & $n_{\rm gmres}$ & $t_{\mathrm{comp}}$
  & $n_{\rm gmres}$ & $t_{\mathrm{comp}}$
  & $n_{\rm gmres}$ & $t_{\mathrm{comp}}$
 & $n_{\rm gmres}$ \\ 
\hline
& 16  & 0.08 & 242   & 0.05 & 41    & 0.06 & 138   & 0.05 & 31  \\  
1& 32  & 0.28 & 1256  & 0.47 & 874   & 0.07 & 281   & 0.06 & 49  \\  
& 64  & 0.58 & 2786  & -- & -- & 0.36 & 1501  & 0.10 & 83  \\  
& 128 & 1.27 & 4485  & -- & -- & 0.86 & 3184  & 0.30 & 169 \\  
\hline
& 16  & 0.04 & 203   & 0.04 & 37    & 0.03 & 126   & 0.05 & 27  \\  
2& 32  & 0.14 & 555   & 0.07 & 62    & 0.07 & 240   & 0.06 & 39  \\  
& 64  & 1.18 & 4449  & 0.15 & 129   & 0.28 & 1273  & 0.12 & 77  \\  
& 128 & 1.48 & 4450  & 0.61 & 292   & 1.13 & 3266  & 0.42 & 169 \\  
\hline
& 16  & 0.06 & 240   & 0.05 & 37    & 0.04 & 150   & 0.05 & 32  \\  
3& 32  & 0.18 & 707   & 0.07 & 68    & 0.08 & 256   & 0.07 & 49  \\  
& 64  & 0.53 & 2545  & 0.16 & 137   & 0.29 & 1244  & 0.14 & 84  \\  
& 128 & 1.65 & 4559  & 0.91 & 412   & 1.25 & 3188  & 0.55 & 195 \\  
\hline
& 16  & 0.05 & 233   & 0.05 & 41    & 0.04 & 149   & 0.05 & 30  \\  
4 & 32  & 0.14 & 729   & 0.08 & 69    & 0.07 & 267   & 0.08 & 48  \\  
& 64  & 0.59 & 2603  & 0.26 & 156   & 0.35 & 1319  & 0.18 & 88  \\  
& 128 & 2.73 & 6415  & 2.44 & 982   & 1.88 & 3779  & 1.22 & 348 \\  
\hline
\end{tabular}
\caption{Computational time $t_{\mathrm{comp}}$
 (in seconds) and total number of GMRES iterations $n_{\rm gmres}$  as a function of $k$ and $h$ for BJ, BJ-PP(10), ASM, and ASM-PP(10) preconditioners. The reported GMRES iteration counts correspond to the total number of linear iterations accumulated over all Newton iterations required to converge the steady-state problem.}
\label{tab1}
\end{table}

The growth of $n_{\mathrm{gmres}}$ with mesh refinement is significantly mitigated by the ASM-based methods, which exhibit only a moderate increase in iteration counts as $1/h$ doubles. This behavior contrasts sharply with the BJ preconditioner, for which the iteration counts increase considerably with the number of elements. Similarly, increasing the polynomial degree $k$ leads to only a modest rise in iteration counts for ASM-PP(10), indicating good robustness with respect to both $h$- and $p$-refinement. Hence, ASM-PP(10) delivers scalable performance across a broad range of discretization parameters.

Table \ref{tab2} decomposes the total computational time into four components — assembly time $t_{\mathrm{ass}}$, matrix–vector product time $t_{\mathrm{mv}}$, preconditioner application time $t_{\mathrm{prec}}$, and orthogonalization time $t_{\mathrm{orth}}$ — for both the ASM and ASM-PP(10) preconditioners. These quantities provide a detailed view of the relative computational costs associated with different stages of the GMRES solver and how they evolve with increasing polynomial degree $k$ and mesh refinement $1/h$. For both preconditioners, the assembly cost $t_{\mathrm{ass}}$ is relatively small and grows moderately with $k$ and $1/h$, reflecting the local elementwise nature of the HDG formulation. The matrix–vector product time $t_{\mathrm{mv}}$ and the preconditioner time $t_{\mathrm{prec}}$ dominate the total runtime as the mesh is refined, with both quantities approximately doubling when $1/h$ is halved. This behavior indicates good computational scaling and confirms that the matrix operations and domain-decomposition preconditioners scale nearly linearly with the number of degrees of freedom.

\begin{table}[htbp]
\centering
\begin{tabular}{|c|c|cccc|cccc|}
\hline
  \multicolumn{2}{|c|}{} & \multicolumn{4}{c|}{ASM} & \multicolumn{4}{c|}{ASM-PP(10)}   \\
\hline
$k$ & $1/h$ & $t_{\mathrm{ass}}$
  & $t_{\rm mv}$ & $t_{\mathrm{prec}}$
  & $t_{\rm orth}$ & $t_{\mathrm{ass}}$
  & $t_{\rm mv}$ & $t_{\mathrm{prec}}$
 & $t_{\rm orth}$ \\ 
\hline
& 16  &  0.004  &  0.003  &  0.003  &  0.048  &  0.004  &  0.001  &  0.014  &  0.008  \\  
1& 32  &  0.006  &  0.007  &  0.007  &  0.054  &  0.006  &  0.001  &  0.021  &  0.054  \\  
& 64  &  0.008  &  0.043  &  0.042  &  0.256  &  0.009  &  0.003  &  0.054  &  0.016  \\  
& 128  &  0.03  &  0.160  &  0.155  &  0.493  &  0.03  &  0.009  &  0.18  &  0.052  \\  
\hline 
& 16  &  0.005  &  0.003  &  0.003  &  0.017  &  0.005  &  0.001  &  0.013  &  0.008  \\  
2 & 32  &  0.008  &  0.006  &  0.006  &  0.047  &  0.008  &  0.001  &  0.018  &  0.013  \\  
& 64  &  0.014  &  0.045  &  0.044  &  0.17  &  0.014  &  0.003  &  0.061  &  0.015  \\  
& 128  &  0.053  &  0.242  &  0.243  &  0.571  &  0.053  &  0.013  &  0.285  &  0.031  \\  
\hline
& 16  &  0.009  &  0.004  &  0.004  &  0.02  &  0.007  &  0.001  &  0.015  &  0.009  \\  
3 & 32  &  0.01  &  0.007  &  0.006  &  0.055  &  0.01  &  0.001  &  0.022  &  0.014  \\  
 & 64  &  0.021  &  0.044  &  0.043  &  0.17  &  0.023  &  0.003  &  0.069  &  0.016  \\  
& 128  &  0.082  &  0.262  &  0.260  &  0.615  &  0.082  &  0.017  &  0.362  &  0.042  \\  
\hline 
& 16  &  0.009  &  0.004  &  0.004  &  0.02  &  0.008  &  0.001  &  0.014  &  0.008  \\  
4 & 32  &  0.015  &  0.008  &  0.007  &  0.038  &  0.015  &  0.002  &  0.03  &  0.012  \\  
&  64  &  0.041  &  0.059  &  0.057  &  0.184  &  0.043  &  0.004  &  0.09  &  0.016  \\  
& 128  &  0.158  &  0.450  &  0.450  &  0.789  &  0.158  &  0.042  &  0.882  &  0.076  \\  
\hline
\end{tabular}
\caption{Assembly time $t_{\mathrm{ass}}$, matrix-vector product time $t_{\mathrm{mv}}$, preconditioner application time $t_{\mathrm{prec}}$, and operationalization time $t_{\mathrm{orth}}$
 (in seconds) as a function of $k$ and $h$ for ASM, and ASM-PP(10) preconditioners.}
\label{tab2}
\end{table}

The orthogonalization time $t_{\mathrm{orth}}$, which accounts for the Gram–Schmidt process within GMRES, exhibits a more pronounced increase with mesh refinement, particularly for the standard ASM preconditioner. This trend arises from the growth in Krylov subspace dimension as more GMRES iterations are required for convergence. In contrast, ASM-PP(10) consistently achieves much smaller $t_{\mathrm{orth}}$ values, owing to its significantly reduced iteration counts and improved spectral conditioning of the preconditioned system. Comparing the two approaches, ASM-PP(10) generally incurs a modest increase in preconditioner cost $t_{\mathrm{prec}}$ due to the ten polynomial applications per GMRES iteration but compensates for this through substantial savings in matrix–vector and orthogonalization times. As a result, the overall computational balance shifts in favor of ASM-PP(10) for finer meshes and higher-order discretizations, where its superior convergence behavior translates into lower total runtime. These results confirm that ASM-PP preconditioners  deliver scalable and efficient performance.

\subsection{Three-Dimensional Poisson Problem}

We consider the three-dimensional Poisson equation defined on a bounded domain  
\begin{equation}
- \nabla^2 u = 0 \quad \text{in } \Omega \equiv \Omega_{\mathrm{box}} \setminus \Omega_{\mathrm{cone}},
\label{eq:poisson3d}
\end{equation}
where $\Omega_{\mathrm{box}} = (-10,10)^3$ represents a cubic box and $\Omega_{\mathrm{cone}}$ denotes a right circular cone of base radius $R = 0.254$ and height $L = 1.27$ located at the origin $(0,0,0)$ and aligned with the $x$-axis. The cone is subtracted from the box, resulting in a computational domain that features a localized geometric singularity at the cone tip.Homogeneous Dirichlet boundary conditions are imposed on the surface of the cone,  
\[
u = 0 \quad \text{on } \partial \Omega_{\mathrm{cone}},
\]
while a Neumann boundary condition is prescribed on the outer box surface,
\[
\nabla u \cdot \mathbf{n} + (1,0,0) \cdot \mathbf{n} = 0 \quad \text{on } \partial \Omega_{\mathrm{box}},
\]
where $\mathbf{n}$ is the unit outward normal vector. The Neumann condition effectively enforces a unit flux in the $x$-direction. An unstructured tetrahedral mesh consisting of 104{,}429 elements is generated using \texttt{Gmsh}~\cite{Geuzaine2009}, with mesh refinement concentrated near the cone tip to resolve the steep solution gradients induced by the geometric singularity. The numerical solution is computed using polynomial approximation orders $k = 1, 2, 3, 4$ to assess the convergence behavior and  performance of the proposed preconditioning techniques. This three-dimensional problem provides a realistic test case with complex geometry and anisotropic mesh resolution.

The numerical results for the three-dimensional Poisson problem are summarized in Tables~\ref{tab3} and~\ref{tab4}. Table~\ref{tab3} reports the total computational time $t_{\mathrm{comp}}$ and the number of GMRES iterations $n_{\mathrm{gmres}}$ for four different preconditioners. As observed, BJ-PP(10) exhibits slow convergence with 2000 GMRES iterations. Results for the plain BJ and ASM preconditioners are omitted since their performance is quite poor. Increasing the polynomial preconditioner order to $m=20$ dramatically improves convergence for both BJ and ASM methods. For instance, at $k=1$, the GMRES iterations are reduced from 2000 to 89 for BJ and from 2000 to only 37 for ASM. Similar reductions are observed for higher polynomial degrees, confirming that stronger polynomial preconditioning yields a significantly faster convergence. The ASM-based preconditioner consistently outperforms the BJ-based one in terms of both iteration count and wall-clock time. For example, at $k=4$, ASM-PP(20) converges in 53 iterations with a total computation time of 11.48 seconds, while BJ-PP(20) requires 104 iterations and 13.14 seconds. Moreover, the scalability of ASM with respect to the polynomial degree $k$ is better than that of BJ, as reflected by the slower growth in computational time and iteration count with increasing $k$.

Table~\ref{tab4} provides a breakdown of the computational cost for the two best-performing preconditioners, BJ-PP(20) and ASM-PP(20). For both methods, the preconditioner application time $t_{\mathrm{prec}}$ dominates the total cost, accounting for more than 80\% of the overall runtime. The scaling of $t_{\mathrm{prec}}$ with polynomial degree is roughly linear, indicating that the polynomial preconditioner remains computationally efficient even for higher-order discretizations. The matrix–vector and orthogonalization costs remain relatively small and well-controlled, while the assembly time increases moderately with $k$ due to the higher number of degrees of freedom per element. ASM-PP(20) exhibits a smaller preconditioner application time across all polynomial degrees. For instance, at $k=4$, the preconditioner time is 9.86 seconds for ASM versus 11.43 seconds for BJ, reflecting the superior efficiency of the ASM operator. Furthermore, ASM achieves shorter orthogonalization times, consistent with its reduced number of GMRES iterations. 

\begin{table}[h]
\centering
\begin{tabular}{|c|cc|cc|cc|cc|}
\hline
   & \multicolumn{2}{c|}{BJ-PP(10)} & \multicolumn{2}{c|}{BJ-PP(20)} & \multicolumn{2}{c|}{ASM-PP(10)} & \multicolumn{2}{c|}{ASM-PP(20)}   \\
\hline
$k$ & $t_{\mathrm{comp}}$
  & $n_{\rm gmres}$ & $t_{\mathrm{comp}}$
  & $n_{\rm gmres}$ & $t_{\mathrm{comp}}$
  & $n_{\rm gmres}$ & $t_{\mathrm{comp}}$
 & $n_{\rm gmres}$ \\ 
\hline
 1  &  18.947  &  2000  &  1.727  &  89  &  20.5  &  2000  &  0.844  &  37  \\  
 2  &  30.949  &  2000  &  3.203  &  101  &  6.359  &  297  &  1.964  &  43  \\  
 3  &  64.007  &  2000  &  6.806  &  103  &  9.879  &  193  &  5.09  &  48  \\  
 4  &  121.056  &  2000  &  13.137  &  104  &  29.8  &  297  &  11.483  &  53  \\  
\hline
\end{tabular}
\caption{Computational time $t_{\mathrm{comp}}$
 (in seconds) and total number of GMRES iterations $n_{\rm gmres}$  as a function of $k$ for BJ-PP(10), BJ-PP(20), ASM-PP(10), and ASM-PP(20) preconditioners for the 3D Poisson problem.}
\label{tab3}
\end{table}

\begin{table}[h]
\centering
\begin{tabular}{|c|cccc|cccc|}
\hline
   & \multicolumn{4}{c|}{BJ-PP(20)} & \multicolumn{4}{c|}{ASM-PP(20)}   \\
\hline
$k$ & $t_{\mathrm{ass}}$
  & $t_{\rm mv}$ & $t_{\mathrm{prec}}$
  & $t_{\rm orth}$ & $t_{\mathrm{ass}}$
  & $t_{\rm mv}$ & $t_{\mathrm{prec}}$
 & $t_{\rm orth}$ \\ 
\hline
1  &  0.037  &  0.04  &  1.536  &  0.047  &  0.04  &  0.017  &  0.713  &  0.022  \\  
 2  &  0.073  &  0.095  &  2.807  &  0.095  &  0.09  &  0.041  &  1.693  &  0.043  \\  
 3  &  0.187  &  0.234  &  5.973  &  0.145  &  0.245  &  0.11  &  4.45  &  0.073  \\  
 4  &  0.531  &  0.48  &  11.429  &  0.209  &  0.675  &  0.247  &  9.855  &  0.108  \\ 
\hline
\end{tabular}
\caption{Assembly time $t_{\mathrm{ass}}$, matrix-vector product time $t_{\mathrm{mv}}$, preconditioner application time $t_{\mathrm{prec}}$, and operationalization time $t_{\mathrm{orth}}$
 (in seconds) as a function of $k$ for BJ-PP(20), and ASM-PP(20) preconditioners for the 3D Poisson problem.}
\label{tab4}
\end{table}

\subsection{Three-Dimensional Linear Elasticity Beam}

\begin{figure}[htbp]
	\centering
 \includegraphics[width=0.7\textwidth]{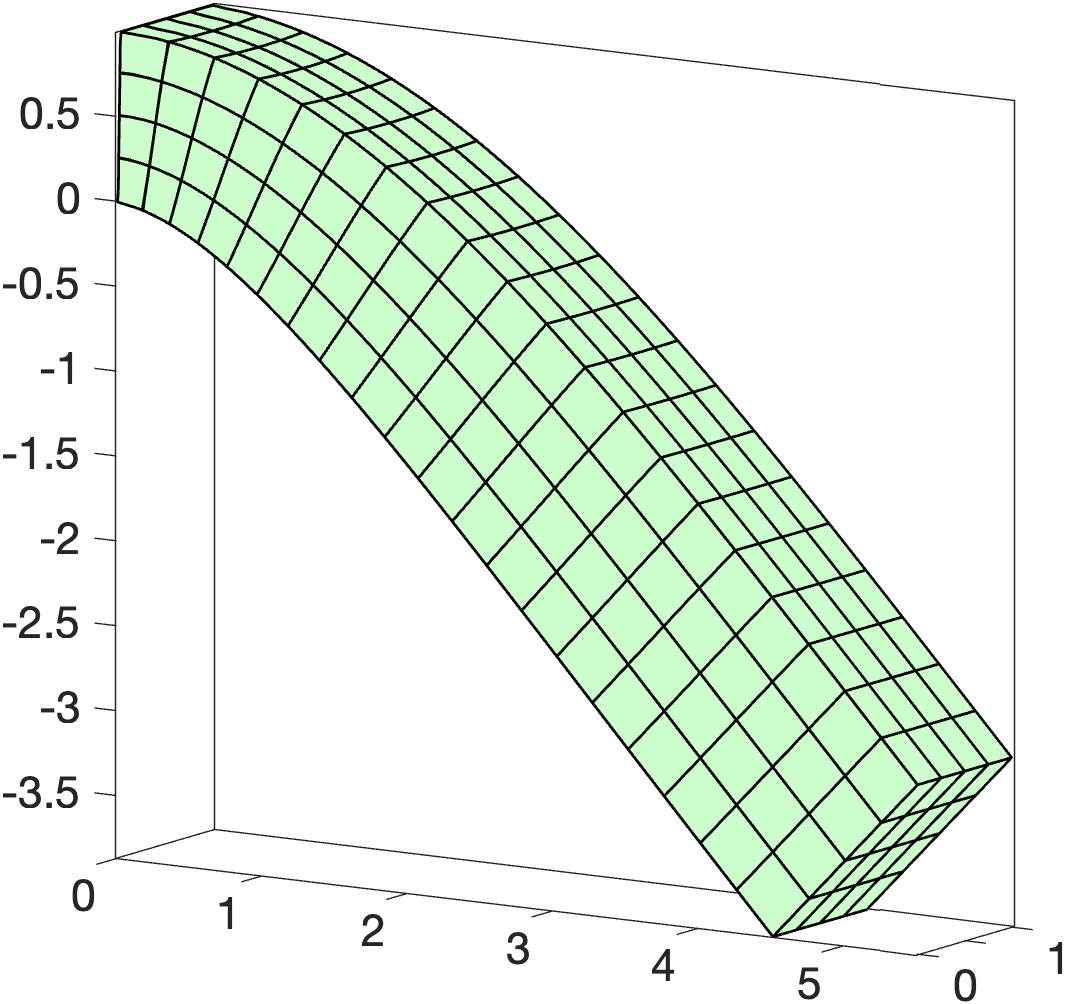}
 \caption{The numerical solution of the 3D linear elasticity beam problem for $k=3$ and $N_E = 320$.}
	\label{fig1}
\end{figure}

We next consider a three-dimensional linear elasticity problem defined on a cantilever beam geometry. The computational domain $\Omega = (0, 5) \times (0,1) \times (0,1)$ represents a rectangular beam of unit thickness, fixed at one end and subjected to a uniform traction load on the opposite face at $x=5$. Homogeneous Dirichlet boundary conditions are imposed on the clamped face, while Neumann conditions are applied to the loaded surface to induce bending deformation. The material is assumed to be isotropic with unit Lamé parameters $\lambda = \mu = 1$, and the displacement field $\mathbf{u} = (u_x, u_y, u_z)$ satisfies the equilibrium equations
\begin{equation}
-\nabla \cdot \boldsymbol{\sigma} = \mathbf{0}, \qquad \boldsymbol{\sigma} = \lambda (\nabla \cdot \mathbf{u}) \mathbf{I} + 2 \mu \,  \boldsymbol{\varepsilon}(\mathbf{u}), \qquad \mbox{in } \Omega,
\end{equation}
where $\boldsymbol{\varepsilon}(\mathbf{u}) = (\nabla \mathbf{u} + \nabla \mathbf{u}^T)/2$ is the strain tensor and $\boldsymbol{\sigma}$ is the Cauchy stress tensor. This problem serves as a benchmark for evaluating solver performance on coupled, vector-valued PDE systems. Several structured hexahedral meshes are employed, consisting of $10 \times 2 \times 2$, $20 \times 4 \times 4$, $30 \times 6 \times 6$, and $40 \times 8 \times 8$ elements. Polynomial approximation orders $k = 1, 2, 3, 4$ are considered to assess the convergence and computational performance of the proposed preconditioning strategies. Figure~\ref{fig1} shows the numerical displacement field for the case $k=3$ and $N_E = 320$ elements, illustrating the smooth bending deformation of the beam under the applied load.

Table~\ref{tab5} presents  $t_{\mathrm{comp}}$ and $n_{\mathrm{gmres}}$ for various preconditioners as a function of polynomial degree $k$ for $N_E = 320$. The unpreconditioned BJ-PP and ASM solvers exhibit extremely slow convergence, requiring thousands of iterations to reach the desired tolerance. In contrast, the polynomially preconditioned versions BJ-PP(10) and ASM-PP(10) achieve dramatic reductions in iteration counts and runtime across all polynomial orders. For instance, at $k=2$, the number of GMRES iterations decreases from 16,000 for BJ-PP to only 141 for BJ-PP(10), and from 3,000 for ASM to merely 24 for ASM-PP(10). The ASM-PP(10) preconditioner consistently delivers the best performance, converging in fewer than 30 iterations for all polynomial degrees with the lowest wall-clock time. Table~\ref{tab6} examines scalability with respect to the number of elements $N_E$. As expected, the computational time increases with mesh refinement, but the iteration counts remain nearly constant for the ASM-PP(10) preconditioner, confirming its robustness with respect to mesh size. The iteration count and runtime of BJ-PP(10) are consistently higher than those of ASM-PP(10). Both BJ-PP and ASM without polynomial preconditioning fail to converge within the iteration limit for larger meshes, highlighting the necessity of polynomial enhancement in three-dimensional elasticity applications. This demonstrates the strong synergy between the additive Schwarz method and polynomial preconditioning in accelerating convergence for elasticity systems.

\begin{table}[htbp]
\centering
\begin{tabular}{|c|cc|cc|cc|cc|}
\hline
   & \multicolumn{2}{c|}{BJ} & \multicolumn{2}{c|}{BJ-PP(10)} & \multicolumn{2}{c|}{ASM} & \multicolumn{2}{c|}{ASM-PP(10)}   \\
\hline
$k$ & $t_{\mathrm{comp}}$
  & $n_{\rm gmres}$ & $t_{\mathrm{comp}}$
  & $n_{\rm gmres}$ & $t_{\mathrm{comp}}$
  & $n_{\rm gmres}$ & $t_{\mathrm{comp}}$
 & $n_{\rm gmres}$ \\ 
\hline
1  &  2.733  &  14000  &  0.082  &  82  &  0.174  &  628  &  0.036  &  19  \\  
 2  &  4.917  &  16000  &  0.277  &  141  &  1.197  &  3000  &  0.096  &  24  \\  
 3  &  5.268  &  9000  &  0.735  &  154  &  2.627  &  3000  &  0.306  &  28  \\  
 4  &  20.277  &  15000  &  2.084  &  162  &  6.241  &  3000  &  0.972  &  29  \\    
\hline
\end{tabular}
\caption{Computational time $t_{\mathrm{comp}}$
 (in seconds) and total number of GMRES iterations $n_{\rm gmres}$  as a function of $k$ for BJ, BJ-PP(10), ASM, and ASM-PP(10) preconditioners for the 3D linear elasticity problem.}
\label{tab5}
\end{table}

\begin{table}[htbp]
\centering
\begin{tabular}{|c|cc|cc|cc|cc|}
\hline
   & \multicolumn{2}{c|}{BJ} & \multicolumn{2}{c|}{BJ-PP(10)} & \multicolumn{2}{c|}{ASM} & \multicolumn{2}{c|}{ASM-PP(10)}   \\
\hline
$N_E$ & $t_{\mathrm{comp}}$
  & $n_{\rm gmres}$ & $t_{\mathrm{comp}}$
  & $n_{\rm gmres}$ & $t_{\mathrm{comp}}$
  & $n_{\rm gmres}$ & $t_{\mathrm{comp}}$
 & $n_{\rm gmres}$ \\ 
\hline
40  &  0.931  &  3000  &  0.063  &  26  &  0.073  &  99  &  0.06  &  14  \\  
 320  &  5.268  &  9000  &  0.735  &  154  &  2.627  &  3000  &  0.306  &  28  \\  
 1080  &  28.636  &  18000  &  6.015  &  498  &  35.519  &  14000  &  1.213  &  42  \\  
 2560  &  --  &  --  &  --  &  --  &  50.698  &  9000  &  3.248  &  54 \\
\hline
\end{tabular}
\caption{Computational time $t_{\mathrm{comp}}$
 (in seconds) and total number of GMRES iterations $n_{\rm gmres}$  as a function of $N_E$ for BJ, BJ-PP(10), ASM, and ASM-PP(10) preconditioners for the 3D linear elasticity problem.}
\label{tab6}
\end{table}

A more detailed breakdown of computational costs is given in Tables~\ref{tab7} and~\ref{tab8}. The dominant contribution to the total runtime arises from the preconditioner application phase, which scales approximately linearly with both polynomial order and problem size. The assembly and matrix–vector operations account for a smaller fraction of the total cost and exhibit smooth scaling behavior. The ASM-PP(10) method achieves consistently lower preconditioner and orthogonalization times compared to BJ-PP(10), reflecting its superior efficiency and convergence properties. The results show that the ASM-PP(10) preconditioner provides a robust and scalable solution strategy for elasticity problems. It maintains rapid convergence with increasing polynomial degree and mesh resolution.

\begin{table}[htbp]
\centering
\begin{tabular}{|c|cccc|cccc|}
\hline
   & \multicolumn{4}{c|}{BJ-PP(10)} & \multicolumn{4}{c|}{ASM-PP(10)}   \\
\hline
$k$ & $t_{\mathrm{ass}}$
  & $t_{\rm mv}$ & $t_{\mathrm{prec}}$
  & $t_{\rm orth}$ & $t_{\mathrm{ass}}$
  & $t_{\rm mv}$ & $t_{\mathrm{prec}}$
 & $t_{\rm orth}$ \\ 
\hline
 1  &  0.005  &  0.003  &  0.046  &  0.021  &  0.006  &  0.001  &  0.016  &  0.008  \\  
 2  &  0.025  &  0.015  &  0.187  &  0.032  &  0.022  &  0.003  &  0.055  &  0.009  \\  
 3  &  0.133  &  0.045  &  0.507  &  0.034  &  0.105  &  0.008  &  0.168  &  0.01  \\  
 4  &  0.769  &  0.109  &  1.132  &  0.04  &  0.526  &  0.02  &  0.386  &  0.011  \\ 
\hline
\end{tabular}
\caption{Assembly time $t_{\mathrm{ass}}$, matrix-vector product time $t_{\mathrm{mv}}$, preconditioner application time $t_{\mathrm{prec}}$, and operationalization time $t_{\mathrm{orth}}$
 (in seconds) as a function of $k$ for BJ-PP(10), and ASM-PP(10) preconditioners for the 3D linear elasticity problem.}
\label{tab7}
\end{table}

\begin{table}[htbp]
\centering
\begin{tabular}{|c|cccc|cccc|}
\hline
   & \multicolumn{4}{c|}{BJ-PP(10)} & \multicolumn{4}{c|}{ASM-PP(10)}   \\
\hline
$N_E$ & $t_{\mathrm{ass}}$
  & $t_{\rm mv}$ & $t_{\mathrm{prec}}$
  & $t_{\rm orth}$ & $t_{\mathrm{ass}}$
  & $t_{\rm mv}$ & $t_{\mathrm{prec}}$
 & $t_{\rm orth}$ \\ 
\hline
40  &  0.016  &  0.002  &  0.019  &  0.015  &  0.03  &  0.001  &  0.016  &  0.008  \\  
 320  &  0.133  &  0.045  &  0.507  &  0.034  &  0.105  &  0.008  &  0.168  &  0.01  \\  
 1080  &  0.44  &  0.471  &  4.855  &  0.16  &  0.335  &  0.04  &  0.778  &  0.02  \\  
 2560  &  --  &  --  &  --  &  --  &  0.741  &  0.118  &  2.268  &  0.031  \\  
\hline
\end{tabular}
\caption{Assembly time $t_{\mathrm{ass}}$, matrix-vector product time $t_{\mathrm{mv}}$, preconditioner application time $t_{\mathrm{prec}}$, and operationalization time $t_{\mathrm{orth}}$
 (in seconds) as a function of $N_E$ for BJ-PP(10), and ASM-PP(10) preconditioners for the 3D linear elasticity problem.}
\label{tab8}
\end{table}

\subsection{Three-Dimensional Nonlinear Elasticity Cube Block}

\begin{figure}[htbp]
	\centering
 \includegraphics[width=0.7\textwidth]{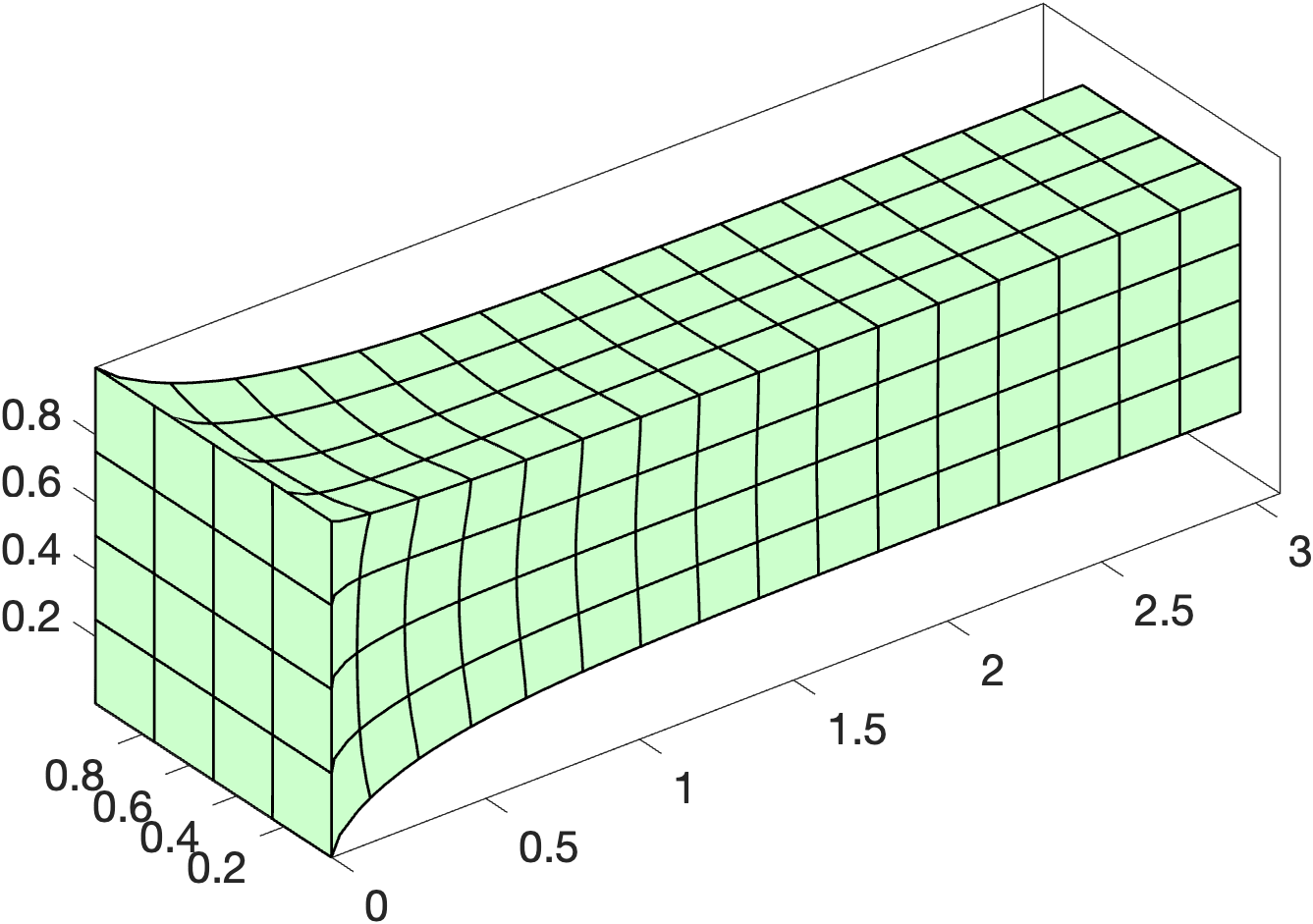}
 \caption{The numerical solution of the 3D nonlinear elasticity problem for $k=3$ and $N_E = 256$.}
	\label{fig3}
\end{figure}

We next examine the performance of the proposed preconditioners for a three-dimensional nonlinear elasticity problem defined on a unit cube discretized with hexahedral elements. The model incorporates large-deformation kinematics and a compressible Neo–Hookean constitutive law. The problem is formulated in terms of the deformation mapping $\bm u(\bm x)$ and its gradient $\bm F = \nabla \bm u$. For a compressible Neo–Hookean material with Lamé parameters $\mu = 1$ and $\lambda = 1$, the deformation gradient $\bm F$ gives rise to the Jacobian $J = \det(\bm F)$ and the right Cauchy–Green deformation tensor $\bm C = \bm F^{T}\bm F$. The strain-energy density is then defined as
\begin{equation}
W(\bm F)
= \frac{\mu}{2}\left(\operatorname{tr}(\bm C) - 3\right)
-		\mu \ln J + \frac{\lambda}{2}(J - 1)^2 .
\end{equation}
The first Piola–Kirchhoff stress tensor follows by differentiating $W$ with respect to $\bm F$. Dirichlet boundary conditions are imposed on one face of the cube and prescribed tractions are applied on the opposite face. Figure~\ref{fig3} shows a representative deformation field computed with polynomial degree $k=3$ and $N_E = 256$ elements.


To study how polynomial degree affects solver performance, Table~\ref{tab9} reports the total computational time and cumulative number of GMRES iterations for $k=1,2,3$, using BJ, BJ-PP(10), ASM, and ASM-PP(10) preconditioners. Several trends emerge.
First, both BJ and BJ-PP(10) exhibit substantially higher iteration counts compared to ASM-based methods, with BJ-PP(10) often requiring more iterations than BJ itself. In particular, for $k=1$, BJ-PP(10) takes more than 6{,}000 iterations, making it the slowest among all tested approaches.
Second, ASM consistently reduces the GMRES iteration counts to roughly one-half or one-third of BJ, and its computational time is correspondingly lower. Finally, ASM-PP(10) delivers the most robust performance across all polynomial degrees: it achieves a dramatic reduction in GMRES iterations—down to only 73–111 iterations—together with the lowest wall-clock time. These results emphasize that the additive Schwarz structure is highly effective for three-dimensional nonlinear elasticity, and that the polynomial enhancement is especially beneficial in this regime.

\begin{table}[htbp]
\centering
\begin{tabular}{|c|cc|cc|cc|cc|}
\hline
   & \multicolumn{2}{c|}{BJ} & \multicolumn{2}{c|}{BJ-PP(10)} & \multicolumn{2}{c|}{ASM} & \multicolumn{2}{c|}{ASM-PP(10)}   \\
\hline
$k$ & $t_{\mathrm{comp}}$
  & $n_{\rm gmres}$ & $t_{\mathrm{comp}}$
  & $n_{\rm gmres}$ & $t_{\mathrm{comp}}$
  & $n_{\rm gmres}$ & $t_{\mathrm{comp}}$
 & $n_{\rm gmres}$ \\ 
\hline
1  &  0.812  &  4603  &  4.015  &  6199  &  0.378  &  1864  &  0.138  &  111  \\  
 2  &  0.769  &  2956  &  3.938  &  2902  &  0.431  &  1101  &  0.267  &  73  \\  
 3  &  1.383  &  2690  &  3.371  &  981  &  1.058  &  1055  &  0.938  &  90  \\  
\hline
\end{tabular}
\caption{Computational time $t_{\mathrm{comp}}$
 (in seconds) and total number of GMRES iterations $n_{\rm gmres}$  as a function of $k$ for BJ, BJ-PP(10), ASM, and ASM-PP(10) preconditioners for the 3D nonlinear elasticity problem.}
\label{tab9}
\end{table}

To quantify scalability with respect to the number of elements, Table~\ref{tab10} examines performance for fixed $k=3$ while increasing $N_E$ from $32$ up to $2048$. The iteration counts for BJ grow rapidly with mesh refinement, reaching over 5{,}000 iterations for the finest mesh. BJ-PP(10) reduces these counts substantially, but its cost per iteration grows enough to make the overall runtime comparable or worse than BJ for large $N_E$. ASM again provides significant improvements, cutting the GMRES counts to about half those of BJ. Most notably, ASM-PP(10) displays the best scalability: GMRES iterations remain below $200$ even for the largest problem, and the total computational time grows much more mildly compared with the other methods. 

\begin{table}[htbp]
\centering
\begin{tabular}{|c|cc|cc|cc|cc|}
\hline
   & \multicolumn{2}{c|}{BJ} & \multicolumn{2}{c|}{BJ-PP(10)} & \multicolumn{2}{c|}{ASM} & \multicolumn{2}{c|}{ASM-PP(10)}   \\
\hline
$N_e$ & $t_{\mathrm{comp}}$
  & $n_{\rm gmres}$ & $t_{\mathrm{comp}}$
  & $n_{\rm gmres}$ & $t_{\mathrm{comp}}$
  & $n_{\rm gmres}$ & $t_{\mathrm{comp}}$
 & $n_{\rm gmres}$ \\ 
\hline
32  &  0.36  &  1235  &  0.235  &  188  &  0.255  &  487  &  0.201  &  50  \\  
 256  &  1.383  &  2690  &  3.371  &  981  &  1.058  &  1055  &  0.938  &  90  \\  
 864  &  5.656  &  4344  &  5.262  &  454  &  5.34  &  2164  &  4.313  &  158  \\  
 2048  &  14.929  &  5519  &  25.435  &  1150  &  15.976  &  3371  &  9.561  &  171  \\  
\hline
\end{tabular}
\caption{Computational time $t_{\mathrm{comp}}$
 (in seconds) and total number of GMRES iterations $n_{\rm gmres}$  as a function of $N_e$ for BJ, BJ-PP(10), ASM, and ASM-PP(10) preconditioners for the 3D nonlinear elasticity problem.}
\label{tab10}
\end{table}

The results show that ASM-PP(10) provides the best overall performance in terms of robustness, iteration counts, and total runtime. The combination of overlapping subdomains and polynomial smoothing yields a highly efficient preconditioner for nonlinear elasticity problems.  Its advantages become increasingly pronounced as the problem size grows, making it particularly effective for both linear and nonlinear elasticity applications.

\subsection{Laminar Subsonic Flow over NACA 0012 Airfoil}

We next consider the two-dimensional laminar, steady, subsonic flow past a NACA~0012 airfoil. The computational domain is constructed using a $C$-type mesh that extends ten chord lengths from the airfoil surface in all directions, ensuring negligible far-field effects. The free-stream Mach number $M_\infty$ and Reynolds number $Re$ are varied to investigate the effect of flow compressibility and viscosity on solver performance. The boundary conditions consist of characteristic inflow and outflow conditions on the far-field boundary and a no-slip adiabatic condition on the airfoil surface. The governing equations are the compressible Navier--Stokes equations in conservative form 
\begin{equation}
\nabla \cdot \bm{F}_{\!c}(\bm{u}) - \nabla \cdot \bm{F}_{\!v}(\bm{u}, \nabla \bm{u}) = \bm{0},
\end{equation}
where $\bm{u} = (\rho, \rho \bm{v}, \rho E)^T$ denotes the vector of conservative variables, and $\bm{F}_{\!c}$ and $\bm{F}_{\!v}$ are the convective and viscous fluxes, respectively. The flow is initialized with the uniform free-stream conditions corresponding to $(M_\infty, Re)$.

\begin{figure}[htbp]
	\centering
 \includegraphics[width=0.49\textwidth]{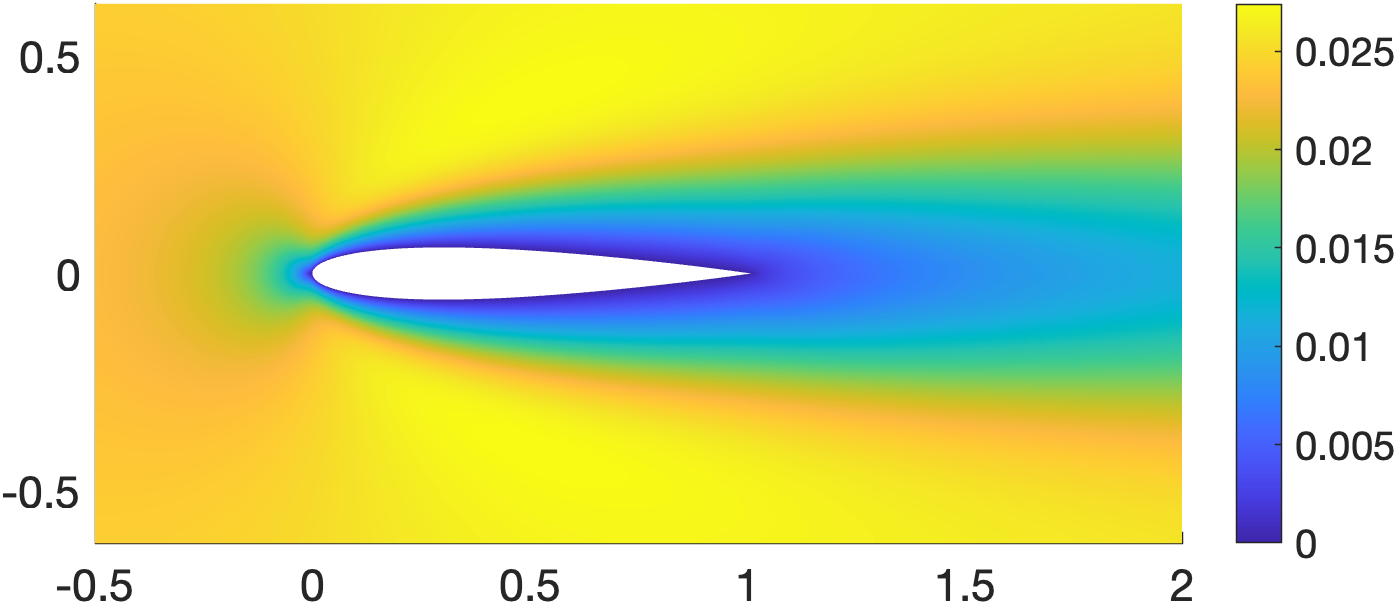}
 \includegraphics[width=0.49\textwidth]{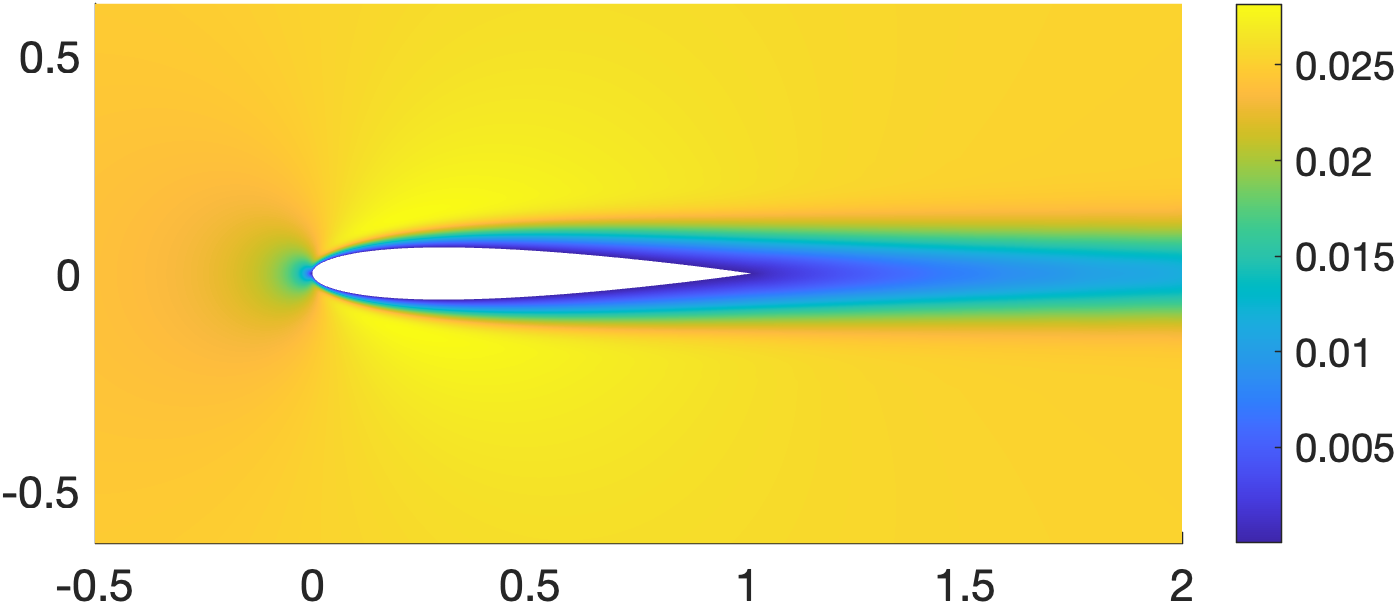}
 \caption{The numerical solution of the subsonic flow over NACA 0012 airfoil for $(M_\infty, Re) = (0.025, 100)$ and $(M_\infty, Re) = (0.025, 1000)$.}
	\label{fig4}
\end{figure}

The simulations are performed using the HDG discretization with polynomial degree $k=3$ on a structured mesh of $1088$ quadrilateral elements. Figure~\ref{fig4} illustrates the computed Mach number contours for $(M_\infty, Re) = (0.025, 100)$ and $(M_\infty, Re) = (0.025, 1000)$. For both cases, the flow remains attached, exhibiting symmetric laminar profiles and smooth pressure recovery along the airfoil surface. The higher-Reynolds-number case produces thinner boundary layers and slightly stronger near-wall velocity gradients, while maintaining the same global flow topology.

Table~\ref{tab11} summarizes the computational time $t_{\mathrm{comp}}$ and GMRES iteration count $n_{\mathrm{gmres}}$ for a range of Mach and Reynolds numbers using the ASM and ASM-PP($m$) preconditioners with $m = 10, 20, 30$. As expected, the convergence behavior is highly sensitive to both $M_\infty$ and $Re$. At low Mach numbers ($M_\infty = 0.025$), the system becomes increasingly ill-conditioned due to the stiffness introduced by the near-incompressible limit, resulting in a large number of GMRES iterations for the unpreconditioned ASM solver (e.g., $n_{\mathrm{gmres}} = 5672$ for $Re = 100$ and $8000$ for $Re = 1000$). The introduction of polynomial preconditioning improves convergence: ASM-PP(10) reduces the iteration count by approximately an order of magnitude, and ASM-PP(20) and ASM-PP(30) further accelerate convergence, reducing the number of iterations by factors of 20–30 relative to the unpreconditioned case. For $Re = 100$ and $M_\infty = 0.025$, the total runtime decreases from 2.23~s for ASM to 1.14~s for ASM-PP(20), representing nearly a twofold speedup. Increasing the polynomial preconditioner order beyond $m=20$ yields diminishing returns: ASM-PP(30) provides slightly fewer GMRES iterations but incurs a marginally higher computational cost per iteration. Overall, ASM-PP(20) achieves the best balance between convergence rate and computational efficiency. 

\begin{table}[htbp]
\centering
\begin{tabular}{|c|c|cc|cc|cc|cc|}
\hline
  \multicolumn{2}{|c|}{} & \multicolumn{2}{c|}{ASM} & \multicolumn{2}{c|}{ASM-PP(10)} & \multicolumn{2}{c|}{ASM-PP(20)} & \multicolumn{2}{c|}{ASM-PP(30)}   \\
\hline
$Re$ & $M_{\infty}$ & $t_{\mathrm{comp}}$
  & $n_{\rm gmres}$ & $t_{\mathrm{comp}}$
  & $n_{\rm gmres}$ & $t_{\mathrm{comp}}$
  & $n_{\rm gmres}$ & $t_{\mathrm{comp}}$
 & $n_{\rm gmres}$ \\ 
\hline
& 0.4  &  0.636  &  1480  &  0.441  &  145  &  0.525  &  86  &  0.593  &  60  \\  
100 & 0.1  &  0.83  &  2019  &  0.551  &  194  &  0.662  &  118  &  0.723  &  82  \\  
& 0.025  &  2.234  &  5672  &  1.74  &  695  &  1.142  &  234  &  1.307  &  177  \\  
\hline 
& 0.4  &  0.601  &  1227  &  0.518  &  151  &  0.589  &  80  &  0.675  &  55  \\  
1000 & 0.1  &  1.165  &  2810  &  0.939  &  344  &  1.059  &  195  &  1.08  &  121  \\  
& 0.025  &  3.059  &  8000  &  3.035  &  1282  &  3.502  &  789  &  3.492  &  521  \\  
\hline
\end{tabular}
\caption{Computational time $t_{\mathrm{comp}}$
 (in seconds) and total number of GMRES iterations $n_{\rm gmres}$ for different values of $(M_\infty, Re)$ for the laminar subsonic flow over NACA 0012 airfoil.}
\label{tab11}
\end{table}

\subsection{Unsteady Euler Vortex}

We consider the two-dimensional unsteady Euler vortex problem  ~\cite{persson2009discontinuous} expressed in conservative form as
\begin{equation}
\frac{\partial \mathbf{u}}{\partial t} + \nabla \cdot \mathbf{F}_{\!c}(\mathbf{u}) = \mathbf{0}, \quad \Omega \equiv (-5,5)^2, 
\end{equation}
where $\mathbf{u} = (\rho, \rho u, \rho v, \rho E)^{T}$ is the vector of conserved variables and $\mathbf{F}_{\!c}$ denotes the inviscid flux tensor. The initial condition consists of a uniform free stream superimposed with an isentropic vortex centered at $(0, 0)$ with prescribed circulation strength and temperature perturbation. The exact analytical solution is periodic in time, allowing direct assessment of the numerical accuracy and stability of the solver without introducing boundary effects. The computational domain is discretized with a quadrilateral mesh of $64 \times 64$ elements, and periodic boundary conditions are imposed on all sides. Time integration is performed using the diagonally implicit Runge-Kutta (DIRK) scheme with fixed time-step sizes $\Delta t = 0.1$, $0.2$, and $0.4$. Polynomial orders $k = 1, 2, 3, 4$ are employed to examine how spatial order and time-step size affects solver convergence and computational efficiency. 

Table~\ref{tab12} reports the assembly time $t_{\mathrm{ass}}$ and GMRES time $t_{\mathrm{gmres}}$ as functions of $\Delta t$ and $k$ for all solvers. As expected, the assembly cost increases approximately linearly with polynomial degree, reflecting the larger number of local degrees of freedom per element, but remains largely insensitive to the time-step size. The GMRES time, on the other hand, varies strongly with both $\Delta t$ and $k$, reflecting changes in stiffness and condition number of the discrete operator. All preconditioners perform comparably well, and polynomial preconditioning provides marginal improvement.The ASM-PP(10) preconditioner achieves the lowest GMRES times across all polynomial degrees.  In contrast, BJ-PP(10) frequently exhibits degraded performance relative to its base counterpart. For example, at $\Delta t = 0.4$ and $k=4$, the GMRES time for BJ-PP(10) rises to 143s—more than twice that of ASM-PP(10), which converges in 65s. These results indicate that the effectiveness of polynomial preconditioning depends strongly on the underlying base preconditioner. While the additive Schwarz formulation benefits from the polynomial acceleration, the block Jacobi variant does not. Overall, the results show that for the unsteady Euler vortex problem, the ASM-PP(10) preconditioner provides the most consistent and efficient performance across a wide range of time-step sizes and polynomial orders.

\begin{table}[htbp]
\centering
\begin{tabular}{|c|c|cc|cc|cc|cc|}
\hline
&   & \multicolumn{2}{c|}{BJ} & \multicolumn{2}{c|}{BJ-PP(10)} & \multicolumn{2}{c|}{ASM} & \multicolumn{2}{c|}{ASM-PP(10)}   \\
\hline
$\Delta t$ & $k$ & $t_{\mathrm{ass}}$
  & $t_{\rm gmres}$ & $t_{\mathrm{ass}}$
  & $t_{\rm gmres}$ & $t_{\mathrm{ass}}$
  & $t_{\rm gmres}$ & $t_{\mathrm{ass}}$
 & $t_{\rm gmres}$ \\ 
\hline
& 1  &  3.73  &  7.29  &  3.74  &  5.62  &  3.89  &  4.83  &  3.77  &  4.98  \\  
0.1 & 2  &  9.49  &  13.39  &  9.68  &  13.66  &  9.55  &  9.51  &  9.76  &  8.83  \\  
& 3  &  20.64  &  17.90  &  20.63  &  49.90  &  20.58  &  13.19  &  20.58  &  12.59  \\  
& 4  &  65.77  &  33.57  &  65.71  &  41.15  &  65.73  &  25.47  &  65.72  &  36.84  \\  
\hline
& 1  &  4.00  &  17.25  &  4.00  &  16.50  &  4.15  &  9.83  &  4.12  &  7.70  \\  
0.2 & 2  &  10.24  &  31.22  &  10.53  &  22.76  &  10.57  &  19.28  &  10.88  &  18.65  \\  
 & 3  &  20.65  &  40.96  &  20.65  &  55.45  &  20.17  &  27.80  &  20.17  &  26.91  \\  
& 4  &  68.39  &  72.19  &  68.41  &  71.47  &  68.57  &  50.49  &  68.56  &  47.23  \\  
\hline 
& 1  &  4.05  &  29.70  &  3.88  &  76.35  &  3.95  &  17.16  &  3.95  &  11.96  \\  
0.4 &  2  &  10.00  &  47.88  &  10.33  &  42.82  &  10.71  &  29.72  &  10.03  &  28.15  \\  
 & 3  &  20.18  &  67.37  &  20.18  &  74.99  &  20.41  &  43.81  &  20.41  &  42.43  \\  
 & 4  &  70.47  &  101.15  &  70.62  &  143.10  &  70.82  &  65.76  &  68.83  &  64.86  \\  
 \hline 
\end{tabular}
\caption{Assembly time $t_{\mathrm{ass}}$ and GMRES time
  $t_{\rm gmres}$  as a function of $\Delta t$ and $k$ for BJ, BJ-PP(10), ASM, and ASM-PP(10) preconditioners for the Euler Vortex problem.}
\label{tab12}
\end{table}

\subsection{Turbulent Flow Past a Flat Plate}

We consider a subsonic turbulent flow over a smooth flat plate at a freestream Mach number $M_\infty = 0.2$ and Reynolds number $Re = 2.1854 \times 10^{6}$~\cite{nguyen07rans,Moro2011a,Moro2017a}. The flow is governed by the compressible Reynolds-Averaged Navier–Stokes (RANS) equations closed with the modified Spalart–Allmaras (SA) turbulence model~\cite{Moro2011a}. The SA model introduces a single transport equation for the modified eddy viscosity variable, which accounts for turbulence production, destruction, and viscous diffusion. The inflow velocity and density correspond to the freestream conditions, while the wall temperature is set to be adiabatic. The mesh consists of 828 quadrilateral elements, refined geometrically toward the wall to accurately capture the steep velocity and turbulent viscosity gradients. A polynomial degree of $k=3$ is employed for both the flow and turbulence variables. Time integration is performed using the implicit backward Euler method with a fixed time-step size $\Delta t$ for 200 time steps.

\begin{figure}[htbp]
	\centering
 \includegraphics[width=0.49\textwidth]{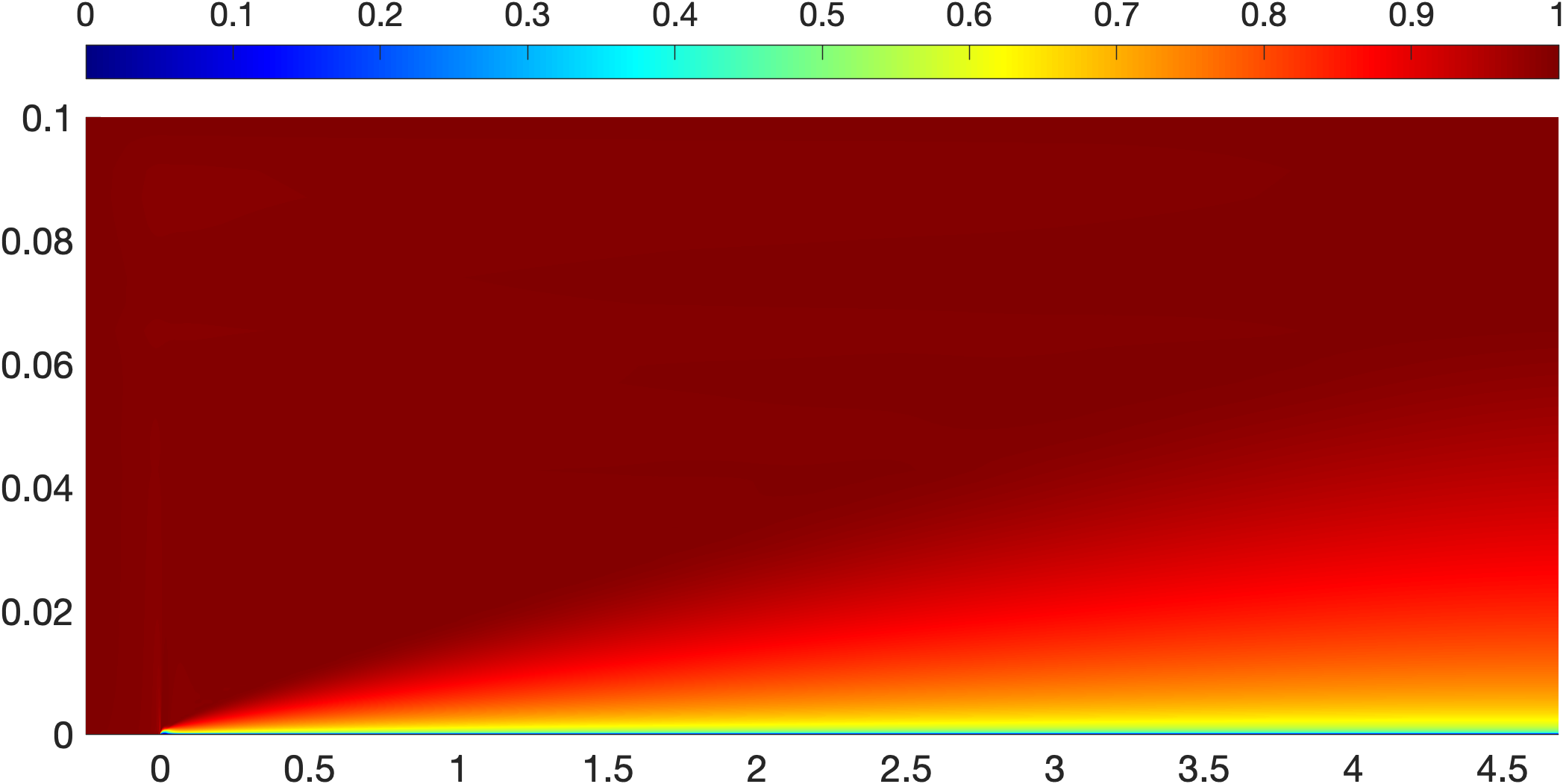}
 \includegraphics[width=0.49\textwidth]{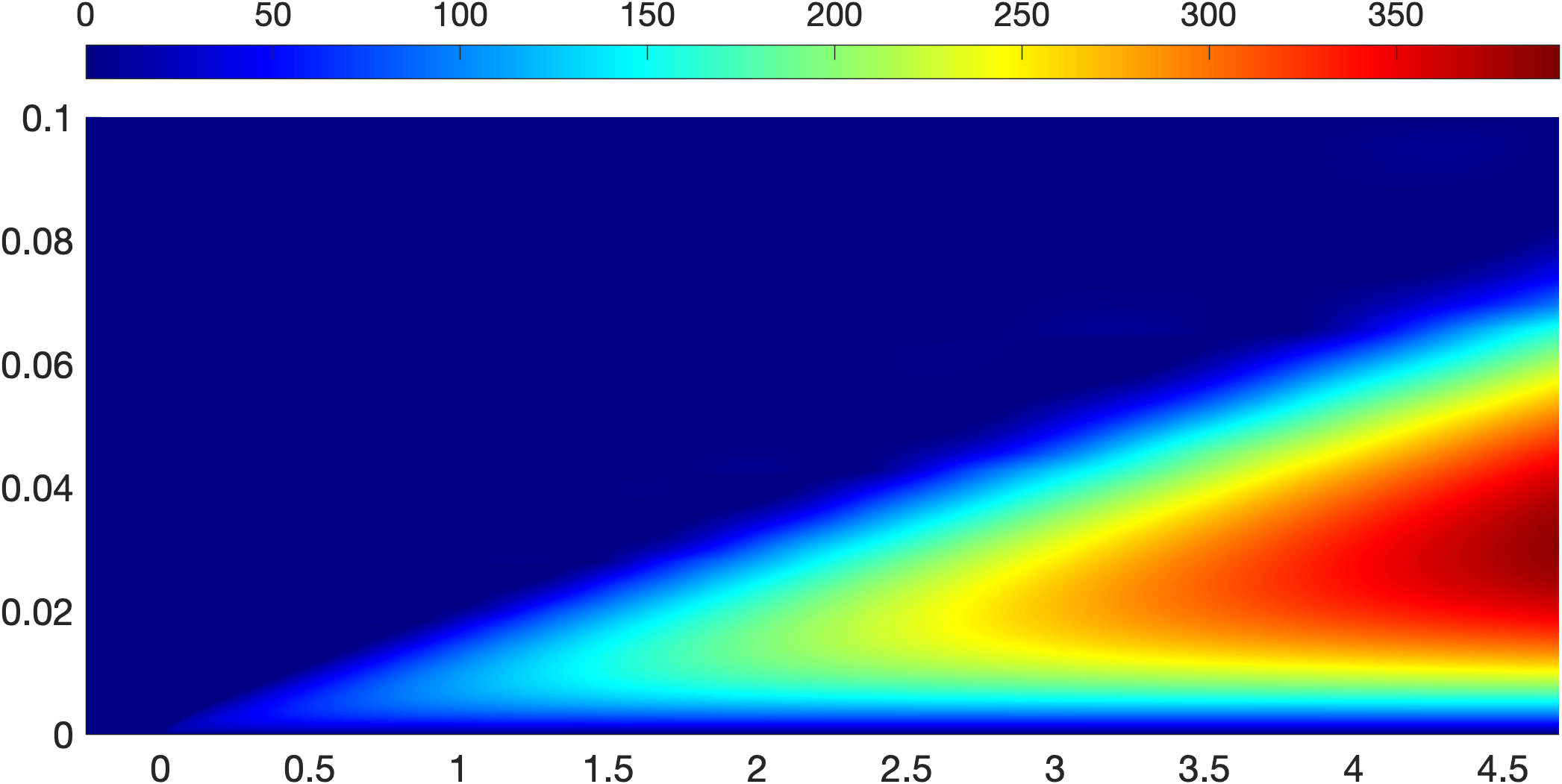}
 \caption{Computed velocity $u/u_{\infty}$ (left) and eddy viscosity ratio $\chi = \tilde{\nu}/\nu$ (right) for the subsonic turbulent flow past flat plate at $(M_\infty, Re) = (0.2, 2.1854 \times 10^6)$.}
	\label{fig5}
\end{figure}

Table~\ref{tab13} summarizes the assembly time $t_{\mathrm{ass}}$ and GMRES time $t_{\mathrm{gmres}}$ as functions of the time-step size $\Delta t$ for different preconditioning strategies.  As expected, the assembly time remains nearly constant across all time steps and preconditioners, since local element operations depend primarily on the polynomial order and mesh size rather than the time-step parameter. In contrast, the GMRES time exhibits a strong dependence on $\Delta t$, reflecting the increased stiffness of the linearized RANS operator as the pseudo–time step grows. For small time steps ($\Delta t \le 4\times10^{-3}$), all preconditioners yield similar convergence behavior with modest GMRES costs. However, as $\Delta t$ increases, the differences in solver robustness become pronounced. Contrary to the laminar case, the results in Table~\ref{tab13} show that polynomial preconditioning does not improve convergence for this turbulent RANS problem. Both ASM-PP(5) and ASM-PP(10) perform worse than the simpler BJ and ASM preconditioners. For instance, at $\Delta t = 1.6\times10^{-2}$, the GMRES time for ASM is only 7.55 seconds, whereas ASM-PP(5) and ASM-PP(10) require 81.19 seconds and 229.83 seconds, respectively. This deterioration becomes even more pronounced at larger time steps. The poor performance of polynomial preconditioning in this case can be attributed to the highly nonlinear and anisotropic nature of the RANS–SA system, which make the polynomial approximation of the inverse operator inaccurate and less effective as a smoother. Flow regimes dominated by turbulence stiffness may require more specialized preconditioning strategies beyond polynomial acceleration.



\begin{table}[htbp]
\centering
\begin{tabular}{|c|cc|cc|cc|cc|}
\hline
   & \multicolumn{2}{c|}{BJ} & \multicolumn{2}{c|}{ASM} & \multicolumn{2}{c|}{ASM-PP(5)} & \multicolumn{2}{c|}{ASM-PP(10)}   \\
\hline
$\Delta t$ & $t_{\mathrm{ass}}$
  & $t_{\rm gmres}$ & $t_{\mathrm{ass}}$
  & $t_{\rm gmres}$ & $t_{\mathrm{ass}}$
  & $t_{\rm gmres}$ & $t_{\mathrm{ass}}$
 & $t_{\rm gmres}$ \\ 
\hline
1e-3  &  2.794  &  5.608  &  2.745  &  2.974  &  2.76  &  2.93  &  2.769  &  4.603  \\  
 2e-3  &  2.763  &  7.287  &  2.759  &  3.745  &  2.815  &  13.55  &  2.806  &  7.623  \\  
 4e-3  &  2.76  &  10.058  &  2.776  &  4.961  &  2.847  &  15.284  &  2.812  &  16.591  \\  
 8e-3  &  2.765  &  13.552  &  2.792  &  6.237  &  2.797  &  16.167  &  2.826  &  21.73  \\  
 1.6e-2  &  2.763  &  18.606  &  2.792  &  7.549  &  2.821  &  81.186  &  2.824  &  229.835  \\  
 3.2e-2  &  2.764  &  23.716  &  2.796  &  9.077  &  2.831  &  187.872  &  2.831  &  347.164  \\  
\hline
\end{tabular}
\caption{Assembly time $t_{\mathrm{ass}}$ and GMRES time
  $t_{\rm gmres}$  as a function of $\Delta t$ for BJ, ASM, ASM-PP(5), and ASM-PP(10) preconditioners for the turbulent flow past a flat plate.}
\label{tab13}
\end{table}

\subsection{Turbulent Flow Past the NACA 0012 Airfoil}

\begin{figure}[htbp]
	\centering
 \includegraphics[width=0.49\textwidth]{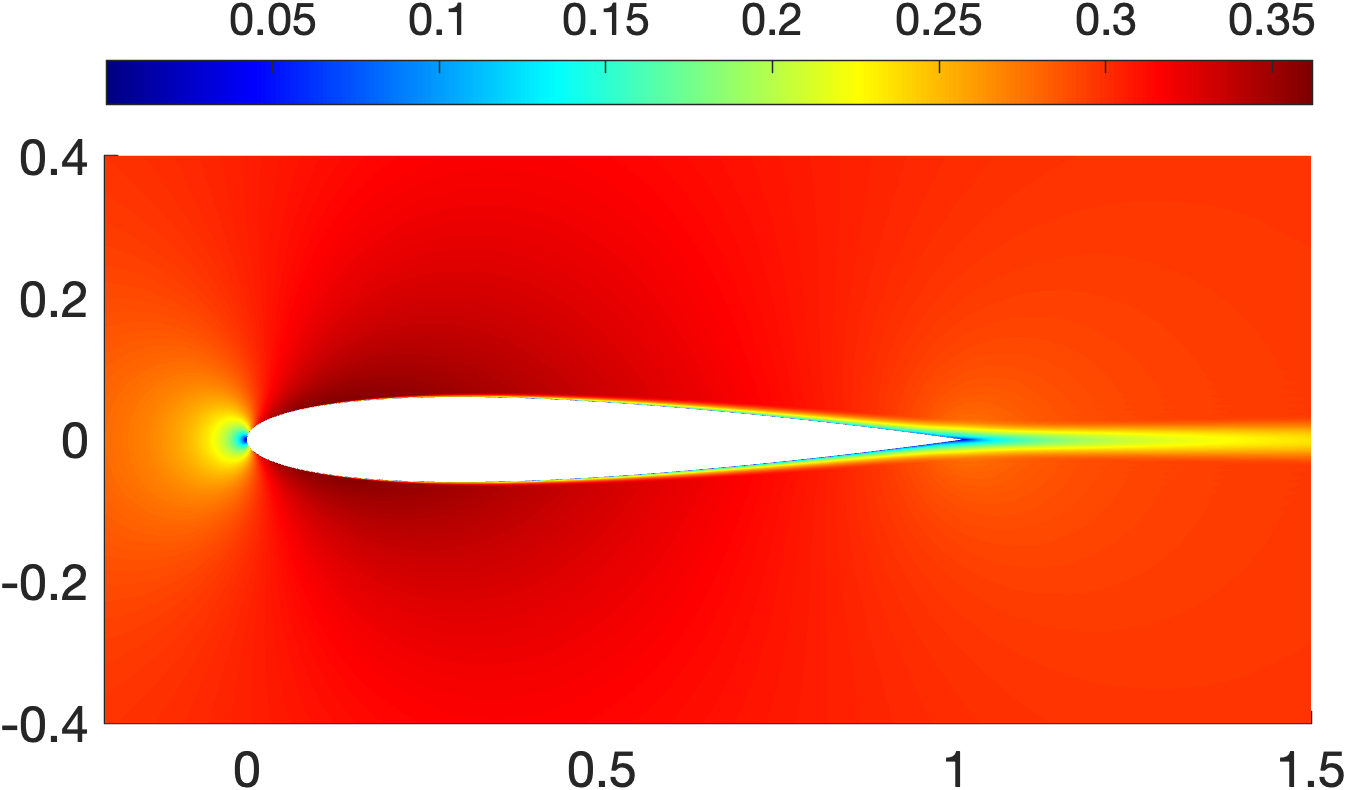}
 \includegraphics[width=0.49\textwidth]{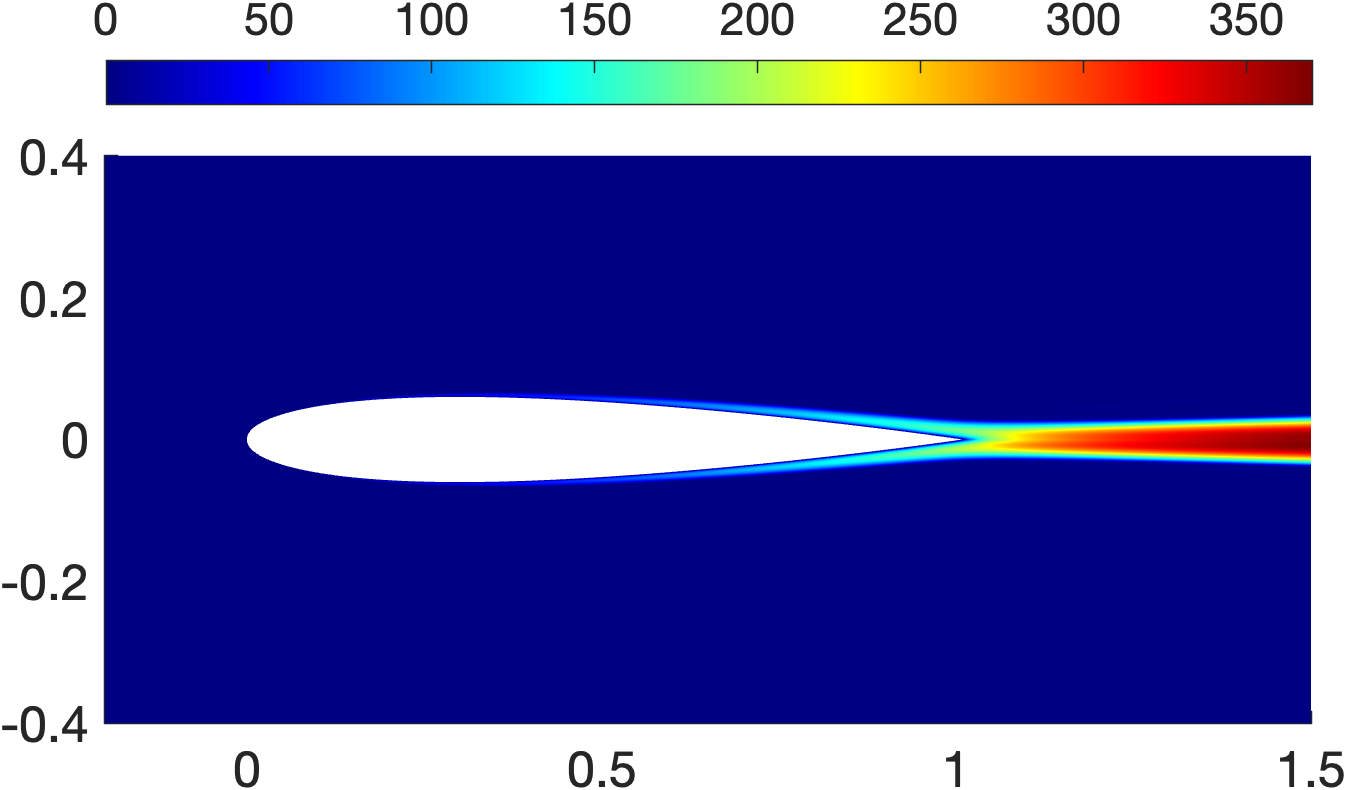}
 \caption{Computed Mach number (left) and eddy viscosity ratio $\chi = \tilde{\nu}/\nu$ (right) for the subsonic turbulent flow past NACA 0012 airfoil at $(M_\infty, Re) = (0.3, 1.85 \times 10^6)$.}
	\label{fig6}
\end{figure}

We present results for a turbulent flow past NACA 0012 airfoil at freestream Mach number $M_\infty = 0.3$, Reynolds number of $1.85 \times 10^{6}$, and zero angle-of-attack~\cite{nguyen07rans,Moro2011a,Moro2017a}. The flow is governed by the RANS equations in conjunction with the modified SA turbulence model~\cite{Moro2011a}. The inflow velocity and density correspond to the freestream conditions, while the wall temperature is set to be adiabatic. We use a single-block, two-dimensional C-grid of 112 x 30 quadrilateral elements. The grid is clustered around the leading edge and the trailing edge to resolve the flow gradients there, and around the airfoil surface to resolve the boundary layer on the airfoil. A polynomial degree of $k=3$ is employed for both the flow and turbulence variables. Time integration is performed using the implicit backward Euler method. Figure~\ref{fig6} displays the computed Mach number field and the eddy-viscosity ratio $\chi = \tilde{\nu}/\nu$, demonstrating a smooth boundary-layer profile and accurate resolution of the wake. 

The performance of the BJ, ASM, ASM-PP(5), and ASM-PP(10) preconditioners is summarized in Table~\ref{tab14}, which reports the average assembly time $t_{\mathrm{ass}}$ and the average GMRES solve time $t_{\rm gmres}$ over a range of time-step sizes $\Delta t$. The assembly cost remains constant across all $\Delta t$ and all preconditioners, reflecting the fact that the residual and Jacobian evaluations are insensitive to the timestep size. In contrast, the GMRES time exhibits a strong dependence on $\Delta t$, with larger timesteps leading to increasingly stiff linear systems and correspondingly higher iteration counts. ASM provides the best performance, reducing GMRES time by a factor of two to five relative to BJ as $\Delta t$ increases. The polynomially enhanced variants show mixed behavior: ASM-PP(5) offers modest improvements at intermediate timesteps but becomes less efficient for very small or very large $\Delta t$, while ASM-PP(10) consistently underperforms ASM. These results indicate that for this RANS-SA problem, the base ASM preconditioner strikes the best balance between robustness and computational efficiency, whereas polynomial preconditioning offers limited benefit.

\begin{table}[htbp]
\centering
\begin{tabular}{|c|cc|cc|cc|cc|}
\hline
   & \multicolumn{2}{c|}{BJ} & \multicolumn{2}{c|}{ASM} & \multicolumn{2}{c|}{ASM-PP(5)} & \multicolumn{2}{c|}{ASM-PP(10)}   \\
\hline
$\Delta t$ & $t_{\mathrm{ass}}$
  & $t_{\rm gmres}$ & $t_{\mathrm{ass}}$
  & $t_{\rm gmres}$ & $t_{\mathrm{ass}}$
  & $t_{\rm gmres}$ & $t_{\mathrm{ass}}$
 & $t_{\rm gmres}$ \\ 
\hline
 1.6e-3  &  48.17  &  9.07  &  48.09  &  6.80  &  48.1  &  11.45  &  48.13  &  17.93  \\  
 3.2e-3  &  48.37  &  20.42  &  48.53  &  13.85  &  48.53  &  20.84  &  48.51  &  26.83  \\  
 6.4e-3  &  48.59  &  58.00  &  48.98  &  31.53  &  48.94  &  36.32  &  48.97  &  60.01  \\  
 1.28e-2  &  48.77  &  131.77  &  49.43  &  43.42  &  49.47  &  51.05  &  49.38  &  109.73  \\  
 2.56e-2  &  48.86  &  285.63  &  49.83  &  55.22  &  49.82  &  90.84  &  49.81  &  127.39  \\  
 5.12e-2  &  49.00  &  400.52  &  50.26  &  56.86  &  50.22  &  123.58  &  50.22  &  160.00  \\
\hline
\end{tabular}
\caption{Average assembly time $t_{\mathrm{ass}}$ and average GMRES time
  $t_{\rm gmres}$  as a function of $\Delta t$ for BJ, ASM, ASM-PP(5), and ASM-PP(10) preconditioners for the turbulent flow past NACA 0012 airfoil.}
\label{tab14}
\end{table}

\subsection{Taylor--Green Vortex}

Finally, we consider the Taylor--Green vortex at $Re = 1600$~\cite{Fernandez2017,VilaPerez2022}, a widely used benchmark for assessing the performance of high-order discretizations and iterative solvers in transitional and weakly turbulent regimes. The problem is solved in the triply periodic domain $\Omega = (0,2\pi)^3$ using a structured mesh of $28^3$ hexahedral elements and polynomial degree $k=2$ for all solution variables. All simulations are performed on a single AMD MI250X GPU, providing a representative test of the efficiency and scalability of the proposed preconditioners in a fully three-dimensional, time-dependent flow setting.

Table~\ref{tab15} reports the average assembly time $t_{\mathrm{ass}}$ and the average GMRES solve time $t_{\rm gmres}$ over a range of time-step sizes $\Delta t$. As expected, the assembly cost is nearly independent of $\Delta t$ and remains comparable across all preconditioners, since the flux evaluations and residual computations do not change with the stiffness of the implicit system. In contrast, the GMRES time is sensitive to the timestep size: as $\Delta t$ increases, the linear systems become more stiff, resulting in longer solve times for all methods. Among the tested preconditioners, ASM consistently delivers the lowest GMRES time across nearly all $\Delta t$, reducing the solve time by a factor of two or more relative to BJ. The polynomially enhanced variants exhibit mixed performance: ASM-PP(5) provides modest gains over BJ but does not outperform ASM, while ASM-PP(10) is significantly more expensive due to the increased cost of applying the polynomial filter. Overall, the results indicate that standard ASM offers the most favorable balance between robustness and computational efficiency for the Taylor--Green vortex at this resolution.

To further understand the behavior of the preconditioners, Table~\ref{tab16} breaks down the GMRES cost into matrix--vector product time ($t_{\mathrm{mv}}$), preconditioner application time ($t_{\mathrm{prec}}$), and orthogonalization time ($t_{\mathrm{orth}}$) as a function of the polynomial degree $k$. Several important trends emerge. First, for all methods, the cost of matrix--vector products increases with $k$, reflecting the higher number of degrees of freedom per element. Second, the ASM preconditioner roughly doubles the cost of its application compared with BJ, but this increase is offset by the substantially reduced iteration count, leading to faster overall solves. Third, polynomial preconditioning (ASM-PP(5)) dramatically increases $t_{\mathrm{prec}}$ with $k$ while offering only limited reductions in $t_{\mathrm{mv}}$ or iteration count. This imbalance leads to significantly higher total costs, explaining its suboptimal performance relative to ASM. Finally, $t_{\mathrm{orth}}$ remains relatively small for all methods, especially on the GPU, indicating that orthogonalization is not a dominant cost in this problem at the selected discretization order.

Together, these results confirm that for moderately high Reynolds numbers and practical polynomial degrees, the classical ASM preconditioner provides the most effective smoothing and convergence acceleration for implicit high-order DG discretizations of the Taylor--Green vortex. Polynomial preconditioning offers limited benefit in this regime, as its additional computational overhead outweighs any improvements in iteration count, particularly in large-scale GPU-accelerated simulations.

\begin{table}[htbp]
\centering
\begin{tabular}{|c|cc|cc|cc|cc|}
\hline
   & \multicolumn{2}{c|}{BJ} & \multicolumn{2}{c|}{ASM} & \multicolumn{2}{c|}{ASM-PP(5)} & \multicolumn{2}{c|}{ASM-PP(10)}   \\
\hline
$\Delta t$ & $t_{\mathrm{ass}}$
  & $t_{\rm gmres}$ & $t_{\mathrm{ass}}$
  & $t_{\rm gmres}$ & $t_{\mathrm{ass}}$
  & $t_{\rm gmres}$ & $t_{\mathrm{ass}}$
 & $t_{\rm gmres}$ \\ 
\hline
1e-3  &  3.31  &  0.04  &  3.54  &  0.04  &  3.56  &  0.12  &  3.56  &  0.23  \\  
 2e-3  &  3.41  &  0.37  &  3.64  &  0.21  &  3.65  &  0.42  &  3.65  &  0.41  \\  
 4e-3  &  3.44  &  0.48  &  3.76  &  0.31  &  3.77  &  0.44  &  3.77  &  0.66  \\  
 8e-3  &  3.50  &  0.30  &  3.87  &  0.32  &  3.89  &  0.48  &  3.89  &  0.59  \\  
 1.6e-2  &  3.53  &  0.45  &  3.97  &  0.36  &  3.98  &  0.50  &  3.98  &  0.69  \\  
 3.2e-2  &  3.55  &  0.78  &  4.09  &  0.56  &  3.91  &  0.86  &  3.91  &  1.04  \\  
 \hline
\end{tabular}
\caption{Average assembly time $t_{\mathrm{ass}}$ and average GMRES time
  $t_{\rm gmres}$  as a function of $\Delta t$ for BJ, ASM, ASM-PP(5), and ASM-PP(10) preconditioners for the Taylor-Green vortex problem.}
\label{tab15}
\end{table}

\begin{table}[htbp]
\centering
\begin{tabular}{|c|ccc|ccc|ccc|}
\hline
   & \multicolumn{3}{c|}{BJ} & \multicolumn{3}{c|}{ASM} & \multicolumn{3}{c|}{ASM-PP(5)}   \\
\hline
$k$  & $t_{\rm mv}$ & $t_{\mathrm{prec}}$
  & $t_{\rm orth}$ & $t_{\rm mv}$ & $t_{\mathrm{prec}}$
 & $t_{\rm orth}$ & $t_{\rm mv}$ & $t_{\mathrm{prec}}$
 & $t_{\rm orth}$ \\ 
\hline
 1  &  0.04  &  0.00  &  0.00  &  0.02  &  0.02  &  0.00  &  0.02  &  0.10  &  0.00  \\  
 2  &  0.29  &  0.03  &  0.05  &  0.10  &  0.10  &  0.01  &  0.04  &  0.38  &  0.00  \\  
 4  &  0.37  &  0.04  &  0.07  &  0.15  &  0.15  &  0.01  &  0.04  &  0.40  &  0.00  \\  
 8  &  0.24  &  0.03  &  0.03  &  0.15  &  0.16  &  0.01  &  0.04  &  0.44  &  0.00  \\  
 16  &  0.35  &  0.04  &  0.06  &  0.17  &  0.17  &  0.02  &  0.05  &  0.45  &  0.00  \\  
 32  &  0.58  &  0.07  &  0.13  &  0.26  &  0.26  &  0.04  &  0.11  &  0.75  &  0.00  \\  
\hline
\end{tabular}
\caption{Matrix-vector product time $t_{\mathrm{mv}}$, preconditioner application time $t_{\mathrm{prec}}$, and operationalization time $t_{\mathrm{orth}}$ (in seconds) as a function of $k$ for BJ, ASM, and ASM-PP(5) preconditioners for the Taylor-Green vortex problem.}
\label{tab16}
\end{table}

\section{Conclusions}

This work presented an efficient implementation and an extensive evaluation of block Jacobi, additive Schwarz, and polynomially enhanced preconditioners for large-scale HDG discretizations on modern GPU architectures. Through a sequence of test cases we assessed the robustness, efficiency, and scalability of these preconditioners across a broad range of polynomial degrees, mesh resolutions, and problem stiffness levels. Several clear trends emerged. First, the ASM preconditioner consistently provided substantial reductions in GMRES iteration counts relative to BJ, often decreasing the nonlinear or implicit solve time by factors of two to ten depending on the problem configuration. This improvement was particularly pronounced in the three-dimensional nonlinear elasticity and turbulent RANS cases, where the discrete operators exhibit strong coupling and stiffness. In these settings, ASM demonstrated excellent robustness to both mesh refinement and polynomial degree, making it an effective default choice for most applications.

Second, the role of polynomial preconditioning was found to be strongly problem dependent. For the Poisson problem and three-dimensional nonlinear elasticity, the addition of a low-order polynomial filter to ASM (e.g., ASM-PP(10)) significantly enhanced convergence, yielding the best overall performance with iteration counts that were nearly mesh- and degree-independent. However, for turbulent flows, polynomial preconditioners offered limited benefit. Their increased application cost often outweighed any iteration count reductions, and in some regimes they were less efficient than the base ASM method. These results highlight the need to tailor polynomial preconditioning to the spectral properties of the underlying operator rather than applying it uniformly across problem classes.

Third, ASM maintains an advantageous balance between matrix–vector product cost and preconditioner application cost on modern accelerators. Polynomial preconditioning significantly increases the latter, especially at higher polynomial degrees, and therefore must be applied judiciously in performance-critical GPU workflows. Orthogonalization costs were shown to remain relatively small across all problems, indicating that the dominant savings arise from reducing preconditioner or matrix–vector cost rather than modifying the Krylov basis operations. Therefore, ASM is a robust, scalable, and GPU-friendly preconditioning strategy for HDG discretizations across a wide range of PDEs, while polynomial enhancements can provide substantial gains for certain elliptic and nonlinear systems. The combination of HDG's compact stencil, GPU-suitable element-local structure, and ASM's effectiveness makes the approach particularly promising for large-scale simulations.

Future work would explore adaptive polynomial preconditioning strategies, incorporate multilevel or coarse-grid corrections, and extend these techniques to hybrid CPU–GPU and multi-GPU environments. Incomplete LU (ILU)–type preconditioners  are widely used in large-scale CFD solvers, and exploring GPU-efficient ILU factorizations could provide an alternative path to high-performance preconditioning for HDG methods. A particularly promising direction is the integration of ILU  factorizations with the ASM framework. In the present work, ASM employs one-element subdomains to maximize GPU locality and enable fully parallel application. By enriching the subdomain solves with ILU factorizations applied on multi-element subdomains, one can potentially achieve significantly stronger smoothing and improved spectral properties while retaining the parallel structure of ASM. Efficient GPU-oriented ILU-ASM preconditioners would offer a scalable and high-performance preconditioning strategy for HDG discretizations for stiff multiphysics problems.

\section*{Acknowledgements} \label{}
We would like to thank Professors Jaime Peraire and Wesley L. Harris at MIT for fruitful discussions. We gratefully acknowledge  the United States Department of Energy under contract DE-NA0003965 for supporting this work, as well as the MIT SuperCloud and Lincoln Laboratory Supercomputing Center for providing HPC resources.

\bibliographystyle{elsarticle-num} 
\bibliography{library.bib}

\end{document}